\documentclass[12pt]{article}
\usepackage[amsmath]{e-jc}

\usepackage{amssymb,amsfonts,amsmath,stmaryrd,bbm}
\usepackage{latexsym}
\usepackage{verbatim}
\usepackage{stmaryrd}
\usepackage{algorithm}
\usepackage{algpseudocode}
\usepackage{colortbl}

\usepackage[utf8]{inputenc}
\usepackage[T1]{fontenc}

\usepackage{kpfonts}
\usepackage{afterpage}
\usepackage{dsfont}
\usepackage{color}
\usepackage{dsfont}
\usepackage{url}
\usepackage{pifont, tikz, subfigure}
\usetikzlibrary{plotmarks,shapes,arrows,positioning}

\usepackage{microtype}

\newcommand{\LandauO}{\mathcal{O}}

\def\qee{$\hfill{\Box}$}

\newcommand{\precdot}{<\mathrel{\mkern-5mu}\mathrel{\cdot}}

\newcommand{\beq}{\begin{equation}}
  \newcommand{\eeq}{\end{equation}}

\renewcommand{\epsilon}{\varepsilon}

\newcommand{\ns}{\mathbb{N}}

\newcommand{\zs}{\mathbb{Z}}

\newcommand{\qs}{\mathbb{Q}}

\newcommand{\bm}[1]{\mbox{\boldmath \ensuremath{#1}}}

\newcommand{\bD}{\bm D}
\newcommand{\bU}{\bm U}

\newcommand{\cP}{\mathcal P}
\newcommand{\cS}{\mathcal S}

\newcommand{\cW}{\mathcal W}

\newcommand\oeis[1]{\href{https://oeis.org/#1}{#1}}

\DeclareMathOperator{\NInc}{\mathtt{NInc}}
\DeclareMathOperator{\NDec}{\mathtt{NDec}}
\DeclareMathOperator{\Low}{\mathtt{Low}}

\DeclareMathOperator{\Pol}{Pol}

\newcommand{\Qo}{Q^\text{old}}
\newcommand{\Zo}{Z^\text{old}}
\newcommand{\fps}{formal power series}
\newcommand{\gf}{generating function}
\newcommand{\gfs}{generating functions}

\def\emm#1,{{\em #1}}

\newcommand{\lqq}{\trianglelefteq}

\newcommand{\LD}{\mathbb{D}}
\newcommand{\cD}{\mathcal{D}}
\providecommand{\keywords}[1]
{
  \small	
  \textbf{\textit{Keywords---}} #1
}

\dateline{Sep 24, 2024}{Feb 17, 2025}{TBD}


\Copyright{The authors. Released under the CC BY-ND license (International 4.0).}
\title{The ascent lattice on Dyck paths}

\author{Jean-Luc Baril\authornote{1}
\and
Mireille Bousquet-M\'elou\authornote{2}
\and
Sergey Kirgizov\authornote{1}
\and 
Mehdi Naima\authornote{3}
}

\authortext{1}{JLB and SK: LIB, Universit\'e Bourgogne Europe, B.P. 47 870, 21078 Dijon Cedex France, (\email{\{barjl,sergey.kirgizov\}@u-bourgogne.fr}).}

\authortext{2}{MBM: CNRS, LaBRI, Universit\'e de Bordeaux, 351 cours de la Lib\'eration,  F-33405 Talence Cedex, France (\email{bousquet@labri.fr}).}

\authortext{3}{MN: LIP6, Sorbonne Universit\'e, 4 place Jussieu, 75005 Paris, France, (\email{mehdi.naima@lip6.fr}).}

\begin{document}
\maketitle

\begin{abstract} In the Stanley lattice defined on Dyck paths of size $n$,  cover relations are obtained by replacing a valley $DU$ by a peak $UD$. We investigate a greedy version of this lattice, first introduced by Chenevi\`ere,  where cover relations  replace a factor $DU^k D$ by $U^kD^2$. By relating this poset to another poset  recently defined by Nadeau and Tewari, we prove that this  still yields a lattice, which we call the \emm ascent lattice, $\LD_n$.  We then count intervals in $\LD_n$. Their generating function is found to be algebraic of degree~$3$. The proof is based on a recursive decomposition of intervals involving two \emm catalytic, parameters. The solution of the corresponding functional equation is inspired by recent work on the enumeration  of walks confined to a quadrant.  We also consider the order induced in $\LD_{mn}$ on \emm $m$-Dyck paths,, that is, paths in which all ascent lengths are  multiples of $m$, and on \emm mirrored $m$-Dyck paths,, in which all descent lengths are multiples of $m$. The first poset $\LD_{m,n}$ is still a lattice for any~$m$, while the second poset $\LD'_{m,n}$  is only a join semilattice when $m>1$. In both cases, the enumeration of intervals is still described by an equation in two catalytic variables. Interesting connections arise with the \emm sylvester congruence, of Hivert, Novelli and Thibon, and again with walks confined to a quadrant. We combine the latter connection with probabilistic results to give asymptotic estimates of the number of intervals in both $\LD_{m,n}$ and $\LD'_{m,n}$. Their form implies that the \gfs\ of intervals are no longer algebraic, nor even D-finite, when $m>1$.
\end{abstract}
 

\section{Introduction and main results}
\label{sec:intro}

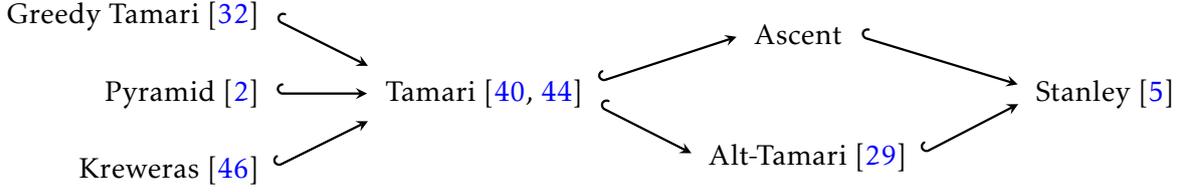
\begin{figure}[htb]
\centering
\scalebox{0.9}{\begin{tikzpicture}[
    mynode/.style={minimum size=4pt},
    inc/.style={line width=0.3mm,
    right hook-stealth,
    shorten <= 3pt,
    shorten >= 3pt}
  ]
  \node[mynode] (TA) at (0, 0)
       {Tamari~\cite{friedman-tamari,HT72}};
  
  \node[mynode, left = 1.4cm of TA] (P)
       {Pyramid~\cite{BKN}};
       
  \node[mynode, above left = 0.38cm and 1.4cm of TA] (D)
       {Greedy Tamari~\cite{aram2022}};
       
  \node[mynode, below left = 0.38cm and 1.4cm of TA] (K)
       {Kreweras~\cite{kreweras-partitions}};
       
  \node[mynode, right = 5.7cm of TA] (S)
       {Stanley~\cite{BeBo07}};
       
  \node[mynode, above right = 0.2cm and 2cm of TA] (A)
       {Ascent};
       
  \node[mynode, below right = 0.2cm and 1.4cm of TA] (alT)
       {Alt-Tamari~\cite{cheneviere1}};

  \draw[inc] (P) to (TA);
  \draw[inc] (D.east) to (TA.north west);
  \draw[inc] (K.east) to (TA.south west);
  \draw[inc] (TA.5) to (A.west);
  \draw[inc] (TA.-5) to (alT.west);
  \draw[inc] (alT.east) to (S.-175);
  \draw[inc] (A.east) to (S.175);

\end{tikzpicture}}
\caption{Subposet-inclusion structure of some orders on Dyck paths.}
\label{figpos}
\end{figure}

In recent years, several orders have been studied on the set of Dyck paths of fixed size~$n$, with a special attention to the number of their intervals. The most natural of these posets is the so-called Stanley lattice (or Dyck lattice~\cite{cheneviere-these}), where a Dyck path is smaller than another one if it lies (weakly) below it~\cite{BeBo07}. Its cover relations are obtained by replacing a \emm valley, by a \emm peak,; in symbols, $DU \rightarrow UD$, where $U$ and $D$ stand for up and down steps, respectively. All other Dyck orders that have been studied are included in this one, in the sense that if a path is less than another one, it lies below it. Let us cite  the Kreweras lattice~\cite{kreweras-partitions,BeBo07}, the Tamari lattice~\cite{friedman-tamari,HT72}, its greedy version due to Dermenjian~\cite{aram2022}, the alt-Tamari lattices of Chenevi\`ere~\cite{cheneviere1}, or the \emm pyramid, lattice introduced by three of the authors of the present paper~\cite{BKN}.
Figure~\ref{figpos} presents the subposet-inclusion structure of these posets.
Interesting connections arise between intervals in these posets and various families of planar maps~\cite{BeBo07,mbm-chapoton,mbm-fusy-preville,ch06,fang-bipartite,fang-bridgeless,fang-trinity,fang-fusy-nadeau}, and, at least conjecturally, with certain quotient rings of polynomials~\cite{bergeron-preville,mbm-chapuy-preville,HaiConj}.

In this paper we consider yet another order
--- in fact, a lattice --- on the set of Dyck paths of size $n$.  This poset was first considered by Chenevi\`ere in his thesis~\cite{cheneviere-these}, following a suggestion of Nadeau.
 The enumeration of its intervals reveals connections with 2-dimensional walks confined to a cone, and, on a more algebraic side, with classes of the \emm sylvester congruence, on words~\cite{HNT,ThibonNovelli}.

 \begin{figure}[htb]
 \centering
\scalebox{0.5}{\begin{tikzpicture}
\draw[line width = 0.8mm, dashed] (0,0) -- (12,0);
\draw[line width = 0.8mm] (0,0) -- (2,2)--(3,1)--(5,3)--(6,2)--(7,3)--(8,4)--(9,3)--(10,2)--(11,1)--(12,0);
\draw[line width = 0.8mm, red] (5,3) -- (6,2)--(8,4)--(9,3);
\filldraw[black] (0,0) circle (4pt);
\filldraw[black] (1,1) circle (4pt);
\filldraw[black] (2,2) circle (4pt);
\filldraw[black] (3,1) circle (4pt);
\filldraw[black] (4,2) circle (4pt);
\filldraw[black] (5,3) circle (4pt);
\filldraw[black] (6,2) circle (4pt);
\filldraw[black] (7,3) circle (4pt);
\filldraw[black] (8,4) circle (4pt);
\filldraw[black] (9,3) circle (4pt);
\filldraw[black] (10,2) circle (4pt);
\filldraw[black] (11,1) circle (4pt);
\filldraw[black] (12,0) circle (4pt);
\end{tikzpicture}}\qquad\scalebox{1.8}{$\precdot$}\qquad
\scalebox{0.5}{\begin{tikzpicture}
\draw[line width = 0.8mm, dashed] (0,0) -- (12,0); 
\draw[line width = 0.8mm] (0,0) -- (2,2)--(3,1)--(5,3);
\draw[line width = 0.8mm] (9,3)--(10,2)--(11,1)--(12,0);
\draw[line width = 0.8mm, red,dashed] (5,3)--(6,2) -- (8,4)--(9,3)--(9,3);
\draw[line width = 0.8mm, red] (5,3) -- (6,4)--(7,5)--(8,4)--(9,3);
\filldraw[black] (0,0) circle (4pt);
\filldraw[black] (1,1) circle (4pt);
\filldraw[black] (2,2) circle (4pt);
\filldraw[black] (3,1) circle (4pt);
\filldraw[black] (4,2) circle (4pt);
\filldraw[black] (5,3) circle (4pt);
\filldraw[black] (6,4) circle (4pt);
\filldraw[black] (7,5) circle (4pt);
\filldraw[black] (8,4) circle (4pt);
\filldraw[black] (9,3) circle (4pt);
\filldraw[black] (10,2) circle (4pt);
\filldraw[black] (11,1) circle (4pt);
\filldraw[black] (12,0) circle (4pt);
\end{tikzpicture}}
\vspace*{3mm}
\caption{A cover relation between two Dyck paths of size $6$.}
\label{figcov}
\end{figure}
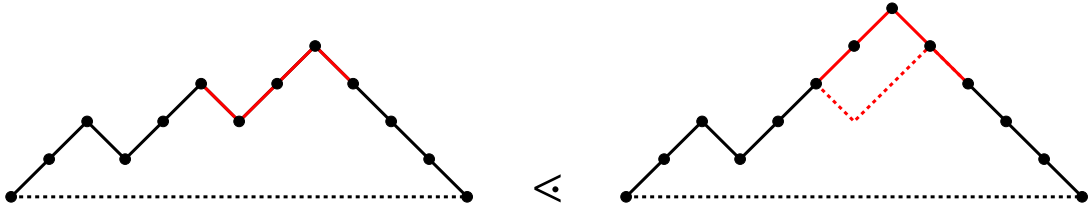

\begin{figure}[htb]
     \centering\includegraphics[scale = 0.45]{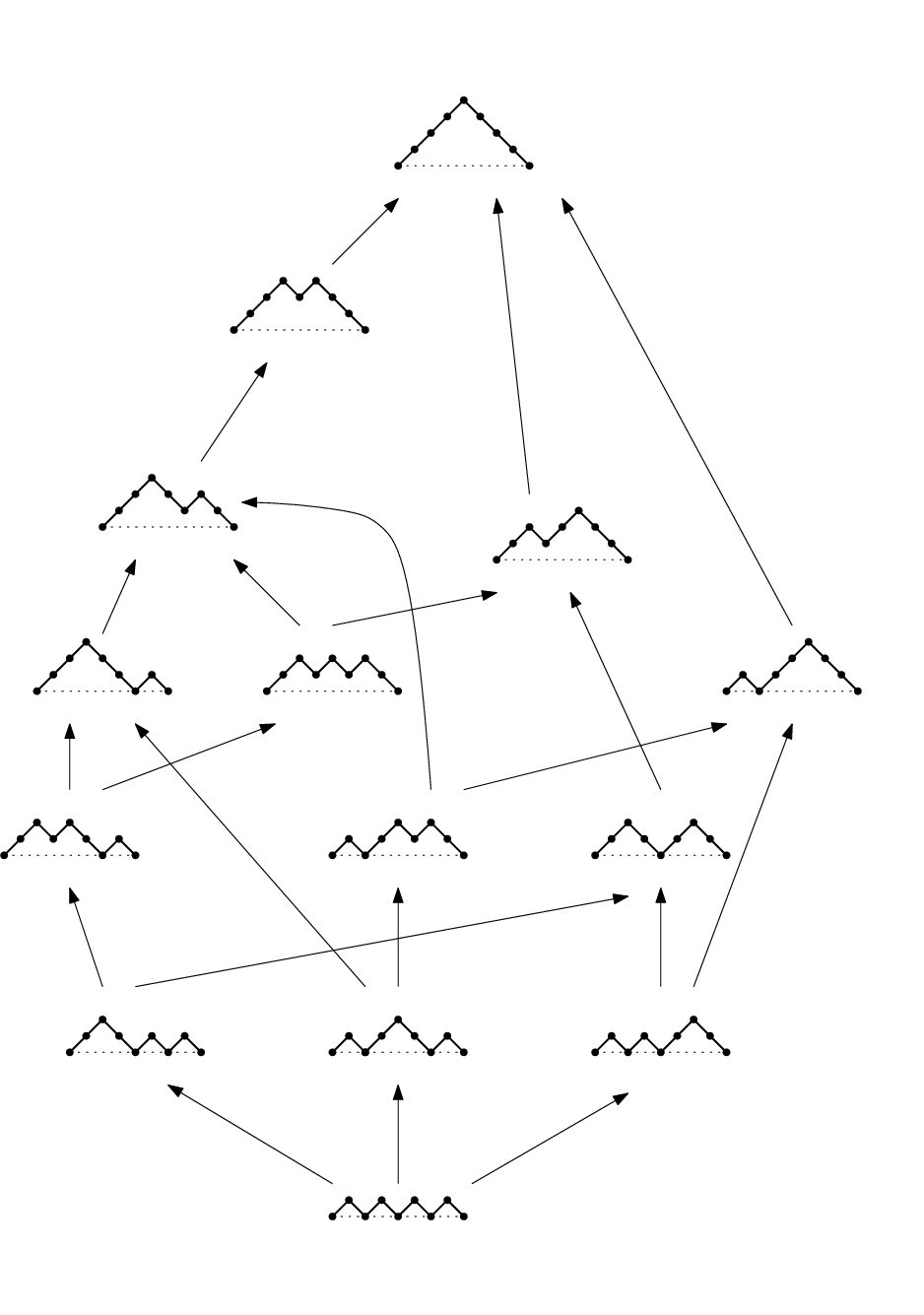}
      \caption{The Hasse diagram of $\LD_4=\LD_{1,4}=\LD'_{1,4}$. }
     \label{fig1}
\end{figure}

The cover relations in this lattice, called the \emm ascent lattice, and denoted by $\LD_n$, are easy to describe. They consist in swapping a down step with the \emm ascent, that follows it, where we call \emm ascent, a (non-empty) maximal sequence of up steps; see Figure~\ref{figcov}.
In symbols, $DU^kD \longrightarrow U^kDD$  for any  $k\geq 1$. Roughly speaking, one applies at once as many cover relations of the Stanley lattice as possible, using always the same down step.  This explains why this poset is called \emm greedy Dyck poset, in~\cite[Def.~7.2.5]{cheneviere-these}. 
Figure~\ref{fig1} shows the whole lattice for $n=4$.
After establishing some properties of the lattice $\LD_n$, we exhibit a recursive construction of its intervals, which can be described by  a \emm generating tree, with two labels. This tree also describes a family of lattice walks confined to the quadrant, in which an \emm infinite family of steps, is allowed. This construction translates, for the associated \gf,  into a linear equation with two additional (or: \emm catalytic,) variables. An important literature has been devoted, in the past 20 years, to the solution of similar-looking equations,  corresponding to quadrant walks in which only \emm finitely many, steps are allowed~\cite{BeBMRa-17,BoKa08,bomi10,DHRS-17,raschel-unified}. We solve our equation by adapting the \emm invariant, approach of~\cite{BeBMRa-17}, and conclude that the \gf\ of intervals in the ascent lattices~$\LD_n$ is an algebraic (cubic) series. The asymptotic behaviour of the corresponding numbers, in $\mu^n n^{-7/2}$, is far less common in enumerative combinatorics than the tree behaviour in $\mu^n n^{-3/2}$, or the (rooted) map behaviour in $\mu^n n^{-5/2}$. This puts ascent intervals in the same universality class as unrooted planar maps~\cite{tutte-polyhedra,richmond-wormald}, or discrete versions of the Brownian motion confined to a wedge of angle $2\pi/5$~\cite{denisov-wachtel,BeBMRa-17}.

\begin{theorem}\label{thm:counting}
  Let $g(n)$ be the number of  intervals in the ascent lattice $\LD_n$. The associated   \gf \ $G 
  :=\sum_{n\ge 1} g(n) t^n$ is
  \[ 
  G=  Z(1- 2Z +2Z^3),
  \]
  where $Z$  is the only \fps\ in $t$ satisfying $Z=t(1+Z)(1 + 2Z)^2$.
  In particular, the series~$G$ is algebraic of degree $3$ over $\qs(t)$.
  
  As $n$ tends to infinity, the number of intervals in $\LD_n$ is equivalent to
  \[
    \kappa\,  \mu^n n^{-7/2},
  \]
with
\[
  \mu= \frac{11+5\sqrt 5}2, \qquad \kappa = \frac 3 {8} \sqrt{\frac{275+123\sqrt 5}{10\, \pi}}.
\]
\end{theorem}
After completing this paper, we learnt that this result had been conjectured by Nadeau and Tewari.

\begin{figure}[htb]
\centering
\scalebox{0.5}{\begin{tikzpicture}
\draw[line width = 0.8mm, dashed] (0,0) -- (12,0);
\draw[line width = 0.8mm] (0,0) -- (2,2)--(3,1)--(5,3)--(6,2)--(7,1)--(9,3)--(10,2)--(11,1)--(12,0);
\filldraw[black] (0,0) circle (4pt);
\filldraw[black] (2,2) circle (4pt);
\filldraw[black] (3,1) circle (4pt);
\filldraw[black] (5,3) circle (4pt);
\filldraw[black] (6,2) circle (4pt);
\filldraw[black] (7,1) circle (4pt);
\filldraw[black] (9,3) circle (4pt);
\filldraw[black] (10,2) circle (4pt);
\filldraw[black] (11,1) circle (4pt);
\filldraw[black] (12,0) circle (4pt);
\end{tikzpicture}}\qquad\qquad
\scalebox{0.5}{\begin{tikzpicture}
\draw[line width = 0.8mm, dashed] (0,0) -- (12,0);
\draw[line width = 0.8mm] (0,0) -- (1,1)--(2,2)--(3,3)--(5,1)--(6,2)--(8,0)--(9,1)--(10,2)--(12,0);
\filldraw[black] (0,0) circle (4pt);
\filldraw[black] (1,1) circle (4pt);
\filldraw[black] (2,2) circle (4pt);
\filldraw[black] (3,3) circle (4pt);
\filldraw[black] (5,1) circle (4pt);
\filldraw[black] (6,2) circle (4pt);
\filldraw[black] (8,0) circle (4pt);
\filldraw[black] (9,1) circle (4pt);
\filldraw[black] (10,2) circle (4pt);
\filldraw[black] (12,0) circle (4pt);
\end{tikzpicture}}
\vspace*{3mm}
\caption{From left to right,
  a 2-Dyck path of size $3$ (element of  $\mathcal{D}_{2,3}$),
   and a mirrored 2-Dyck path of size $3$  (element of  $\mathcal{D}'_{2,3}$).}
\label{figpath}
\end{figure}
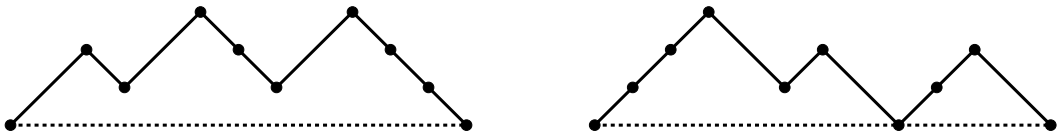

\begin{figure}[htb]
     \centering
     \includegraphics[scale = 0.4]{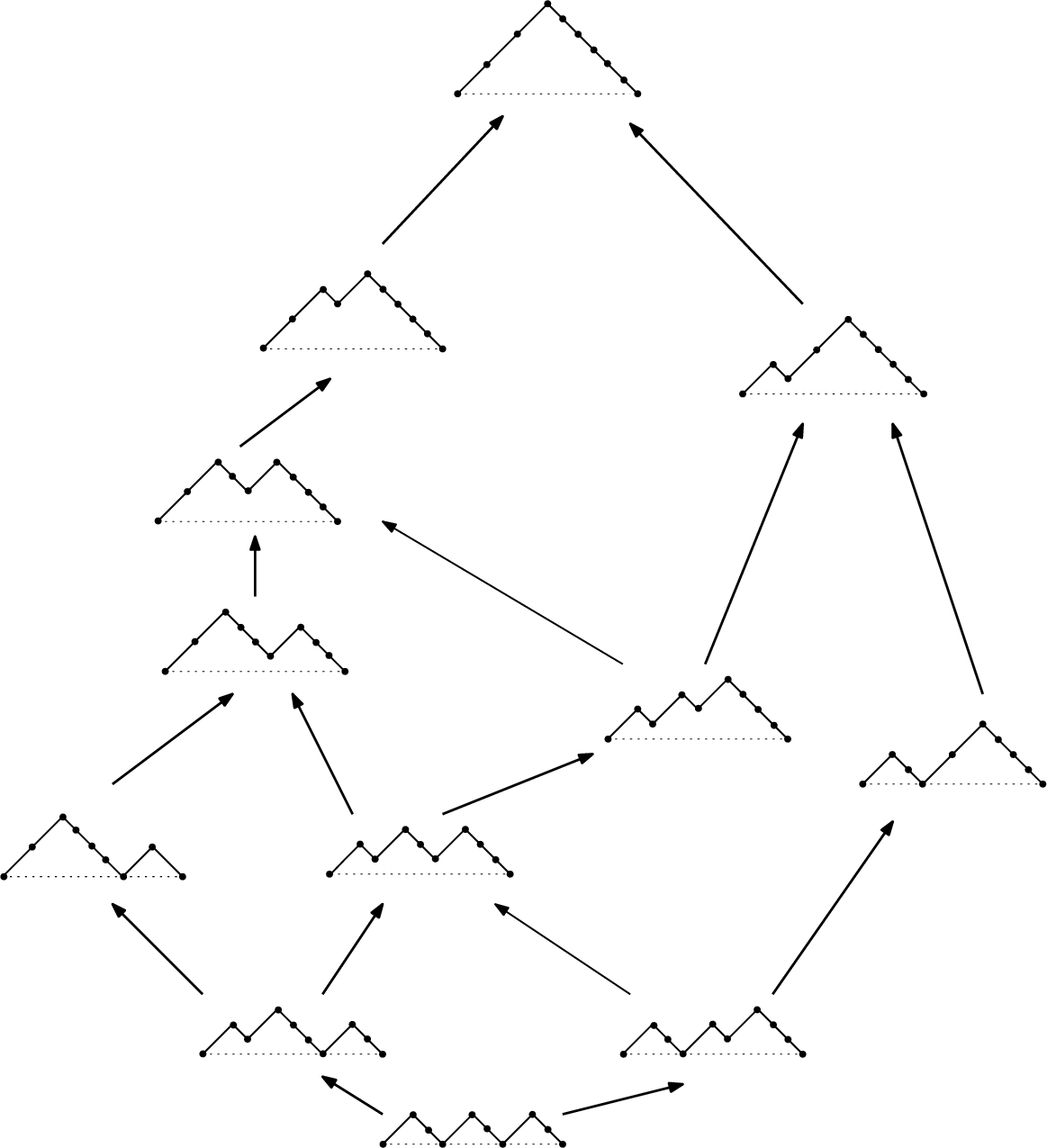}\includegraphics[scale = 0.43]{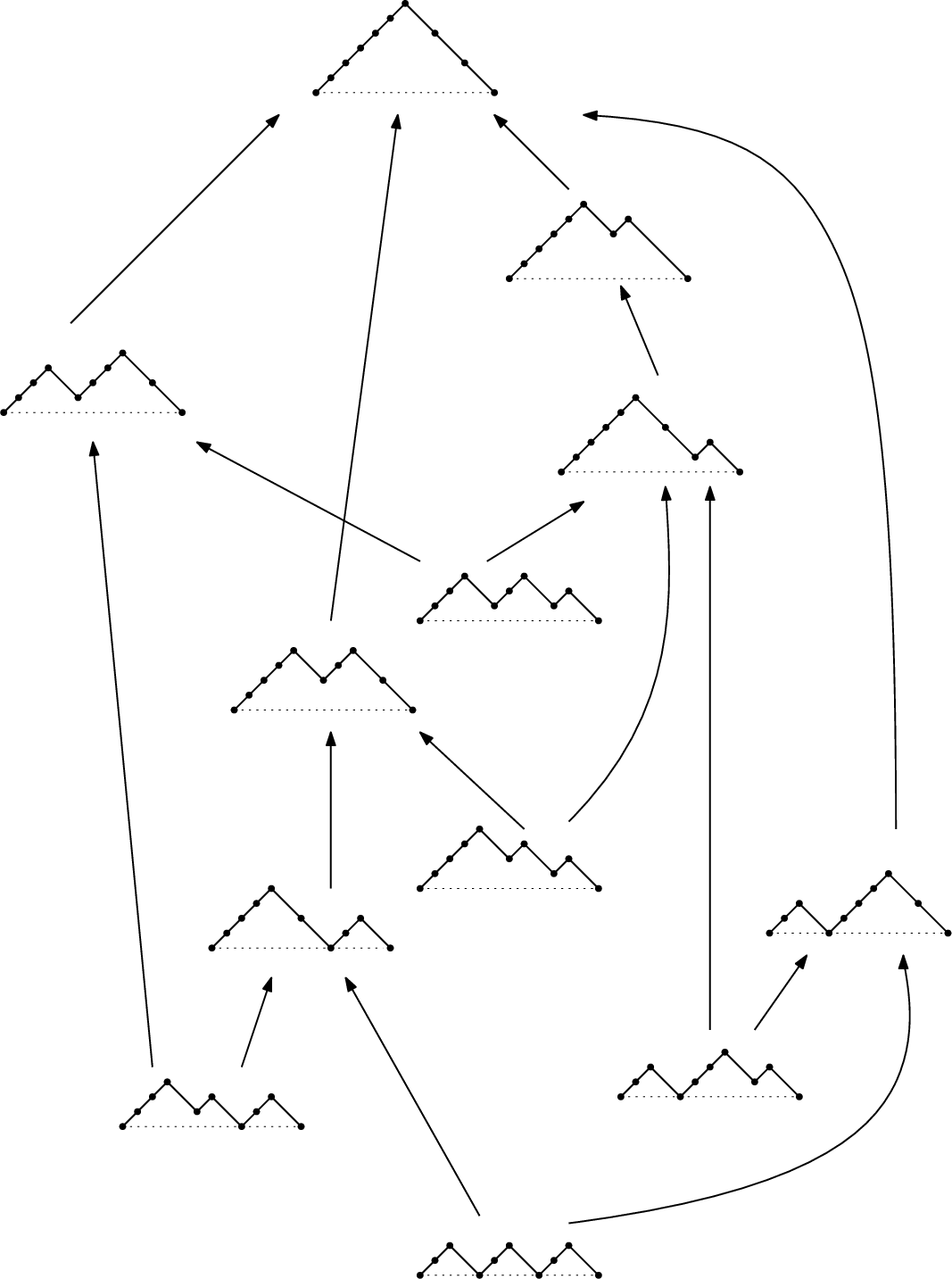}
     \caption{On the left, the Hasse diagram of $\mathbb{D}_{2,3}$, and on the right, the Hasse diagram of~$\mathbb{D}'_{2,3}$. }
     \label{fig2}
 \end{figure}

Inspired by earlier papers on other Dyck lattices~\cite{mbm-chapoton,mbm-fusy-preville,bergeron-preville,mbm-chapuy-preville}, we also consider the family of \emm $m$-Dyck paths, of size $n$, in which the $n$ up steps now have height $m$ instead of $1$ (Figure~\ref{figpath}, left). They can also be seen as Dyck paths of size
$mn$ in which the length of each ascent is a multiple of $m$. They form an upper ideal, and an interval, in $\LD_{mn}$, and thus a new lattice that we denote by $\LD_{m,n}$ (Figure~\ref{fig2}, left). Its cover relations are still described by $D \bm U^k D \longrightarrow \bm U^k DD$, where now $\bm U$, in boldface, stands for a \emm large, up step of height $m$. We extend to this lattice the construction of intervals found for $m=1$, and obtain again a bijection with certain quadrant walks with infinitely many allowed steps. However, in contrast with other  lattices defined on $m$-Dyck paths~\cite{mbm-chapoton,mbm-fusy-preville}, as soon as $m\ge 2$ the associated \gf\ stops being algebraic, or even \emm D-finite, (i.e., solution of a linear differential equation with polynomial coefficients). We prove this by determining the asymptotic behaviour of the number of intervals, which, due to deep number theoretic results~\cite{BoRaSa14}, rules out the possibility of D-finiteness.

\begin{proposition}\label{prop:Dmn}
  Let us fix $m\ge 1$. The number $g_m(n)$ of intervals in the ascent lattice $\LD_{m,n}$ satisfies, as~$n$ tends to infinity,
  \[
     g_m(n) \sim \kappa \mu^n n^{\alpha},
  \]
for some positive constant $\kappa$,   where
  \[
    \mu=\frac{ m \sqrt{m^{2}+4}+m^{2}+2}{2}\cdot \left(\frac{2+\sqrt{m^{2}+4}}{m}\right)^{m}
  \]
  and
  \[
    \alpha=-1 - \pi/\arccos(-c) \qquad \text{with} \qquad c= -\sqrt{\frac{m^{2}+2-\sqrt{m^{2}+4}}{2 m^{2}+6}}.
  \]
  For $m=1$ we have $\alpha=-7/2$, but for $m>1$ the exponent $\alpha$ is irrational. This implies that for $m>1$ the \gf\ of  intervals in $\LD_{m,n}$ is not D-finite.
\end{proposition}

Finally, we also consider the class of \emm mirrored, $m$-Dyck paths of size $n$: now, the up steps have height~$1$ and the $n$ down steps have height $m$ (Figure~\ref{figpath}, right). Alternatively, these paths can be seen as Dyck paths of size
$mn$ in which the length of each descent is a multiple of~$m$. For $m>1$ the order induced by $\LD_{mn}$ on these paths is no longer a lattice, and in particular it has several minimal elements. 
But it is still a join semilattice, which we denote by $\LD'_{m,n}$ (Figure~\ref{fig2}, right). Its cover relations are still described by $\bm D U^k \bm D \longrightarrow U^k \bm {DD}$, where $\bm D$, in boldface, stands for a large down step. Working with these mirrored $m$-Dyck paths is less standard than working with $m$-Dyck paths,  but turns out to be rewarding: when investigating the number of intervals in~$\LD_{m,n}'$, we found them in the On-line Encyclopedia of Integer Sequences~\cite{OnEIS}, as the number of \emm sylvester classes of $m$-parking functions,~\cite{HNT,ThibonNovelli}. We first establish a bijection between these classes and intervals in $\LD'_{m,n}$. Then, we describe a recursive construction of intervals of  $\LD'_{m,n}$, which,  when $m=1$, is \emm not, the same as the earlier one. It gives again a bijection with quadrant walks with infinitely many allowed steps, and a linear equation in two catalytic variables for the \gf. This time we also have an alternative interpretation in terms of quadrant walks with finitely many (weighted) allowed steps. We determine the asymptotic behaviour of the number of intervals, which again rules out the possibility of D-finiteness for $m\ge 2$.

\begin{proposition}\label{prop:Dmnp}
  Let us fix $m\ge 1$. The number $g'_m(n)$ of intervals in the ascent poset $\LD_{m,n}'$ satisfies, as $n$ tends to infinity,
  \[
   g'_m(n) \sim  \kappa \mu^n n^{\alpha},
  \]
for some positive constant $\kappa$,   where
  \[
    \mu=\left(2m+\sqrt{1+4m^2}\right) \left( \frac{1+\sqrt{1+4m^2}}{2m}\right)^{2m}
  \]
  and
  \[
    \alpha=-1 - \pi/\arccos(-c) \qquad \text{with} \qquad c= -\sqrt{\frac{1+2m^2-m\sqrt{1+4m^2}}{2(3m^2+1)}}.
  \]
  For $m=1$ we have $\alpha=-7/2$, but for $m>1$ the exponent $\alpha$ is irrational. This implies that for $m>1$ the \gf\ of ascent intervals in $\LD_{m,n}'$ is not D-finite.
\end{proposition}

\medskip
\noindent
{\bf Outline of the paper.} We begin in Section~\ref{sec:posets} with various definitions, in particular of the ascent posets $\LD_{m,n}$ and $\LD'_{m,n}$. We give a characterization of the ascent order that reveals a connection with another poset defined on nonincreasing sequences, recently introduced by Nadeau and Tewari~\cite{NadeauTewari}. We rely on their results to conclude that $\LD_{m,n}$  is a lattice, and  $\LD_{m,n}'$ a join semilattice.  In Section~\ref{sec:sylvester}, we show how to transform bijectively  intervals of the Nadeau-Tewari poset into words on the alphabet $\zs$ that \emm avoid, two patterns, and are representatives for classes of the
\emm sylvester congruence, introduced by Hivert \emm et al.,~\cite{HNT}. In particular, we exhibit sylvester classes that are in bijection with intervals of  $\LD_{m,n}$ and $\LD'_{m,n}$. In Section~\ref{sec:recursive}, we describe recursive constructions of these  intervals, and convert them into functional equations for the associated \gfs. To write these equations, one needs to record two additional (or: \emm catalytic,) variables. In Section~\ref{sec:m=1} we solve the two equations obtained for $m=1$, and establish Theorem~\ref{thm:counting} and refinements of it. Section~\ref{sec:asympt} is devoted to the asymptotic estimates of Theorems~\ref{prop:Dmn} and~\ref{prop:Dmnp}. We conclude with a few remarks.

\section{Ascent posets}
\label{sec:posets}

\subsection{Definitions}

We begin with precise definitions of the objects and notions discussed more informally in the introduction.

\smallskip\paragraph{\bf Dyck paths.}
Let us first recall that a \emm Dyck path, $P$ is a 2-dimensional path starting at the origin $(0,0)$, consisting of up steps $U=(1,1)$ and  down steps $D=(1,-1)$, that ends on the $x$-axis and never goes strictly below the $x$-axis. The \emm size, of $P$  is the number $n$ of up steps. We denote by $\cD_n$ the set of such paths. An \emm ascent, of a path $P$ is a maximal, non-empty sequence of consecutive up steps. A \emm descent, is defined similarly using down steps. A \emm factor, of $P$ is a non-empty sequence of consecutive steps. A \emm peak, is a factor $UD$, while a \emm valley, is a factor $DU$. The \emm ascent composition, of a Dyck path $P$ is the composition $c(P)=(c_1, \ldots, c_k)$, where the part $c_i>0$ is the length of the $i$th ascent of $P$. Clearly, the $c_i$'s sum to $n$, so that $c(P)$ is a composition of the integer $n$ if $P \in \cD_n$.

For $m\ge 1$, we call \emm $m$-Dyck path, of size $n$ any path of $\cD_{mn}$ in which all ascent lengths are multiples of $m$, and denote the corresponding set by $\cD_{m,n}$. Analogously, we call \emm mirrored $m$-Dyck path, of size $n$ any path of $\cD_{mn}$ in which all descent lengths are multiples of $m$, and denote the corresponding set by $\cD'_{m,n}$. We sometimes consider $m$-Dyck paths of size~$n$ as sequences of \emm large, up steps $\bU=(m,m)$ and (unit) down steps $D$, and analogously for mirrored $m$-Dyck paths, which have large down steps $\bD=(m,-m)$ and unit up steps.
The number of $m$-Dyck paths of size $n$ is the Fuss-Catalan number (see \oeis{A355262}  in the OEIS~\cite{OnEIS}):
\beq\label{fuss}
  C_m(n)=\frac{1}{mn+1}\binom{(m+1)n}{n}.
\eeq

\smallskip\paragraph{\bf Orders on Dyck paths.} There exists on $\cD_n$ a classical order, called the Stanley order (or lattice, in fact), for which $P$ is less than or equal to $Q$ if it lies (weakly) below $Q$. By this, we mean that for any $\ell$, the prefix of $P$ of length $\ell$ contains at most as many up steps as the prefix of length~$\ell$ of $Q$. A path $Q$ covers a path $P$ in this order if and only if $Q$ is obtained by replacing a valley $DU$ of $P$ into a peak $UD$.

In this paper we consider a greedy version of this order, which we call the \emm ascent order,. It is described by its cover relations: we say that $Q$ covers $P$ (denoted  $P\precdot Q$) if there exists in~$P$ a factor $DU^kD$ such that $Q$ is obtained from $P$ by replacing this factor by $U^kDD$. Observe that in this case $P$ lies below $Q$. In particular, this relation is irreflexive and anti-symmetric. Its transitive closure is thus an order relation on $\cD_n$, denoted $\le$, and one easily checks that the cover relations of this order are indeed those described above. We denote by $\LD_n$ the corresponding poset. Since~$\cD_{m,n}$ and $\cD_{m,n}'$ are subsets of $\cD_{mn}$, we can consider  the orders induced by $\le$ on these subsets. We denote by $\LD_{m,n}$ and $\LD'_{m,n}$ the corresponding posets.

\smallskip\paragraph{\bf Formal power series.}
For a ring $R$, we denote by $R[t]$
(resp.~$R[[t]]$) the ring of polynomials (resp. formal power series) 
in $t$ with coefficients in $R$. If $R$ is a field, then $R(t)$ stands for the field of rational functions in~$t$.
This notation is generalized to several variables.
For instance, in Section~\ref{sec:recursive} we consider the \gf\  $G_m(t;x,y)$ of intervals in the ascent lattices $\mathbb{D}_{m,n}$, $n\geq 1$, where $t$ records the size (the number $n$ of up steps), and $x$ (resp. $y$)  the length of the final descent of the smaller (resp. larger) path; this series  belongs to $\qs[x,y ][[t]]$. We often omit in our notation the dependency in $t$ of our series, writing for instance $G_m(x,y)$ instead of $G_m(t;x,y)$.

\subsection{A characterisation of the ascent order}

 Let us first recall that a composition $c$ \emm refines, another composition $c'$ if we can write  $c=(c_{1,1}, c_{1,2}, \ldots, c_{1,j_1}, \ldots, c_{i,1}, c_{i,2}, \ldots, c_{i,j_i})$ and $c'=(c'_1, \ldots, c'_i)$ where $c'_k= c_{k,1}+ c_{k,2}+ \cdots+ c_{k,j_k}$, for $1\le k \le i$.  For instance, the paths of Figure~\ref{figcov} have ascent compositions $c(P)=(2,2,2)$ and $c(Q)=(2,4)$, and we observe that $c(P)$ refines $c(Q)$. Note that refinement is an order relation on compositions of $n$, for each integer $n$. It is isomorphic to the Boolean lattice on $(n-1)$ elements. The cover relations are obtained by merging two consecutive parts.

\begin{proposition}\label{prop:charac}
  In the ascent poset $\LD_n$, we have $P\le Q$ if and only if $P$ lies (weakly) below $Q$ and the ascent composition $c(P)$ refines $c(Q)$.
  
Graphically, the second condition means that for every descent of $Q$, there is a descent of~$P$ lying on the same line (of slope $-1$).\end{proposition}

\begin{proof}
  Let us first prove that if $P\le Q$ in  $\LD_n$, the other two conditions hold. By transitivity, it suffices to prove them if $P\precdot Q$. By definition, this means that one obtains $Q$ from $P$ by replacing a factor $DU^kD$ by $U^k DD$. This shows at once that $Q$ is above $P$. Moreover, $c(Q)$ is either $c(P)$, or is obtained from $c(P)$ by merging two consecutive parts. It follows that $c(P)$ is a refinement of~$c(Q)$.

  \begin{figure}[htb]
 \centering
\hskip -20mm      	\scalebox{0.9}{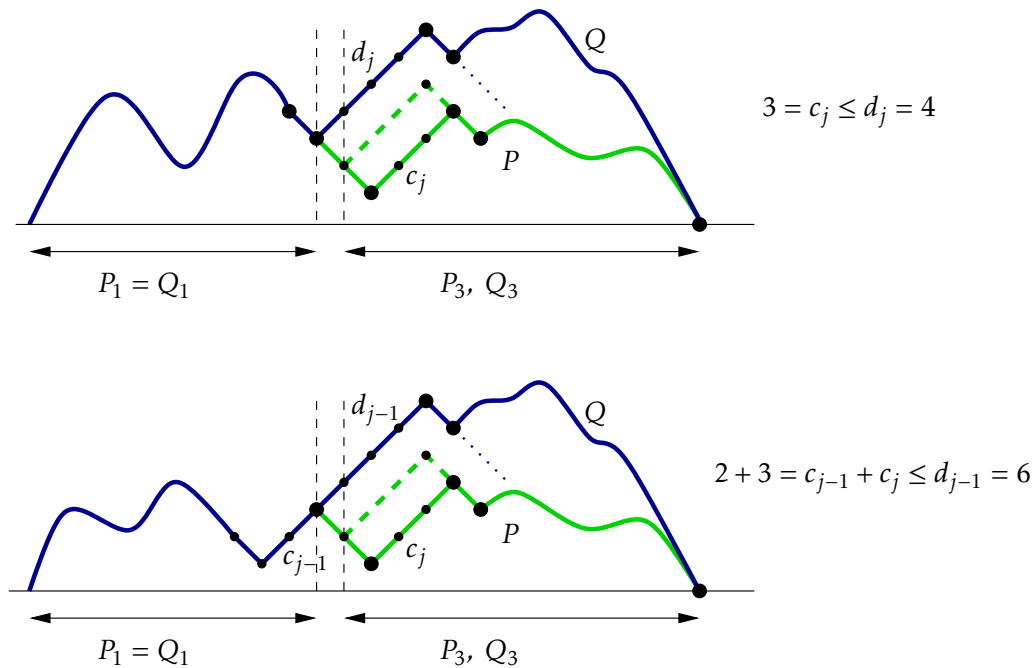}
  \caption{Applying the transformation $DU^kD \rightarrow U^k D D$ at the first valley of~$P$ that is not a valley of $Q$. Top: $P_1$ ends with $D$. Bottom: $P_1$ ends with $U$.
  }
  \label{fig:charac}           
\end{figure}

  Conversely, let us take a pair $(P,Q)$ in $\LD_n$ such that $P$ lies below $Q$ and $c(P)$ refines $c(Q)$. We will argue by induction on the area lying between $P$ and $Q$. If this area is zero, $P=Q$ and there is nothing to prove. Otherwise, let us write $P=P_1 P_2$ and $Q=Q_1 Q_2$ where $P_1=Q_1$ is the longest common prefix of $P$ and $Q$. Then $P_2=DP_3$ and $Q_2=U Q_3$, because $P$ is below~$Q$ (Figure~\ref{fig:charac}).  Note that $P_3$  contains at least one $U$ step. Let $c=(c_1, \ldots, c_k)$ be the ascent composition of $P$, and say that the first ascent of~$P_3$ is the $j$th ascent of $P$, of length $c_j$. If $P_1=Q_1$ ends with a $D$ step (Figure~\ref{fig:charac}, top), then the first up step of $Q_2$ is the start of the $j$th ascent of $Q$, which has length~$d_j$ if $d=(d_1, \ldots, d_\ell)$ is the ascent composition of $Q$. Since~$c$ refines $d$ and the first $(j-1)$ parts of~$c$ and $d$ coincide, we have $c_j \le d_j$. If $P_1=Q_1$ ends with a $U$ step (Figure~\ref{fig:charac}, bottom), then this step belongs to the $(j-1)$th ascents of $P$ and $Q$, which have lengths $c_{j-1}$ in $P$ and $d_{j-1}$ in $Q$. Moreover, $d_{j-1}> c_{j-1}$ because this ascent of $Q$ also includes the first step of $Q_2$ (but definitely not the first step of $P_2$, which is $D$). In this case the first $(j-2)$ parts of $c$ and $d$ coincide, and we have $d_{j-1}\ge c_{j-1}+c_j$ since $c$ refines $d$. In both cases, let $P'$ be obtained by applying the cover relation at the first valley of $P$ that follows $P_1$. The above observations imply that $P'$ lies below $Q$, and that its composition refines $c(Q)$: in the first case, $c(P')=c(P)$, and in the second case, either $c(P')=c(P)$ or $c(P')$ is obtained by merging the parts $c_{j-1}$ and $c_j$ of $c(P)$. Also, the area between $P'$ and $Q$ is less than the area between $P$ and $Q$. By the induction hypothesis, we thus have $P'\le Q$. But $P'$ covers~$P$, hence $P\le Q$.

  The second statement of the proposition is a simple observation.
\end{proof}

It will be worth keeping in mind the following result, established in the second part of the above proof. 

\begin{corollary}\label{cor:first-valley}
  Let $P\lneq Q$.  Applying the cover relation at the first valley of $P$
   that is not a valley of $Q$
  gives a path $P'$ that covers $P$ and satisfies $P'\le Q$.
%
\end{corollary}


\begin{remark}
    The characterization of Proposition~\ref{prop:charac} implies that if $Q$ covers $P$ in the Tamari lattice, then $Q$ is larger than $P$ in the ascent lattice. This inclusion appears in Figure~\ref{figpos}.
\end{remark}


\subsection{The posets $\LD_{m,n}$ and $\LD_{m,n}'$}
We now fix $m\ge 1$. We consider the ascent order on $\cD_{mn}$, and the induced orders on the set $\cD_{m,n}$ of $m$-Dyck paths of size $n$, on the first hand, and on the set $\cD'_{m,n}$ of mirrored $m$-Dyck paths of size $n$, on the other hand. Recall that we sometimes consider these paths as having large up (resp. down) steps, of height $m$, and that these large steps are denoted by $\bU$ and $\bD$, respectively.

\begin{proposition}\label{prop:basic}
  The poset $\LD_{m,n}$ is the interval of $\LD_{mn}$ with minimum element $(U^mD^m)^n$ and maximum element $U^{mn} D^{mn}$. Its cover relations are still given by the transformation $DU^\ell D \rightarrow U^\ell DD$ (where $\ell=mk$ is necessarily a multiple of $m$), or equivalently by $D\bU^kD \rightarrow \bU^k DD$.

  The poset  $\LD_{m,n}'$ has maximum element $U^{mn} D^{mn}$, but several minimal elements if $m\ge 2$. Their number is the Fuss-Catalan number $C_{m-1}(n)$ (see~\eqref{fuss}). The cover relations are  still given by the transformation $\bD U^k\bD \rightarrow U^k \bD\bD$.
\end{proposition}
\begin{proof}
  Recall that the cover relations in $\LD_{mn}$  either merge two consecutive parts of the ascent composition, or leave it unchanged. Hence, if $P\le Q$ in $\LD_{mn}$ and $P\in \cD_{m,n} \subset \cD_{mn}$, then $Q\in \cD_{m,n}$ as well. Thus $\LD_{m,n}$ forms an upper ideal in $\LD_{mn}$. Moreover, Proposition~\ref{prop:charac} implies that the $m$-Dyck path $(U^mD^m)^n$ is smaller than (or equal to) any other $m$-Dyck path. This proves the first statement. The second then follows, because $Q$ covers $P$ in $\LD_{m,n}$ if and only if $Q$ covers~$P$  in~$\LD_{mn}$. 

  Let us now consider the subset $\cD_{m,n}'$ of $\cD_{mn}$. It contains $U^{mn} D^{mn}$, so this path remains the (unique) maximal element in $\LD_{m,n}'$. Let us skip for the moment the results on minimal elements, and focus on cover relations.  If a path $Q$ is obtained from another path $P\in \cD_{m,n}'$ by a transformation $\bD U^k\bD \rightarrow U^k \bD\bD$, then $Q\in \cD'_{m,n}$ and $P\le Q$ by Proposition~\ref{prop:charac}. Moreover, any mirrored $m$-Dyck path that lies above $P$ and below $Q$  is obtained from $P$ by replacing the same factor $\bD U^k \bD$ by $ U^\ell \bD U^{k-\ell}\bD$. By Proposition~\ref{prop:charac}, it is larger than or equal to $P$ in $\LD_{mn}$ only if $\ell=0$ or $\ell=k$. Hence $Q$ covers $P$.
  Conversely, if $Q$ covers $P$ in $\LD_{m,n}'$, then in particular $P\le Q$ in $\LD_{mn}$, and Proposition~\ref{prop:charac} implies that $P$ lies below $Q$ and $c(P)$ refines~$c(Q)$. We  then apply $m$ times the construction of Corollary~\ref{cor:first-valley}: that is, we perform, in the first valley of $P$ that is not a valley of $Q$, a transformation $\bD U^k \bD \rightarrow U^k \bD \bD$ that gives a path $P'$ lying in $[P,Q]$. Since $P'\not = P$ and $P'$ still belongs to $\cD_{m,n}'$, the fact that $Q$ covers $P$ in $\LD'_{m,n}$ implies that $P'=Q$. Hence $Q$ is indeed obtained via the claimed cover relation.

  Now that we have described the cover relations in $\LD'_{m,n}$, let us finally return to its minimal elements. We assume that $m\ge 2$. If $Q$ covers some path $P$, then $Q$ contains a factor $\bD \bD$. Conversely, if $Q= Q_1 U \bD\bD Q_2$, then $Q$ covers $P:=Q_1 \bD U \bD Q_2$. Thus the minimal elements of $\LD'_{m,n}$ are those that contain no factor $\bD \bD$. We claim that these paths are in bijection with mirrored $(m-1)$-Dyck paths of the same size. The bijection simply consists in replacing every factor $U\bD=U D^m$ by $D^{m-1}$. This concludes the proof.
\end{proof}


\subsection{Lattice properties}
\label{sec:lattice}

The characterization of the ascent order in Proposition~\ref{prop:charac} reveals a connection with another order, defined on nonincreasing sequences of integers, recently introduced by Nadeau and Tewari~\cite{NadeauTewari}. This connection, already observed by Chenevi\`ere~\cite[p.~147]{cheneviere-these}, will imply our main structural result.

\begin{definition}{\cite[Def.~5.4]{NadeauTewari}}\label{def:NT}
    Let $n\ge 1$, and $u=(u_1, \ldots, u_n)$, $v=(v_1, \ldots, v_n)$ be two  nonincreasing sequences of integers. Then $u$ is smaller than or equal to $v$ in the Nadeau-Tewari poset $\mathcal P_n$, denoted $u \lqq v$, if $u$ is smaller than or equal to $v$ componentwise, and every descent of $v$ is a descent of $u$. In symbols:
    \begin{itemize}
    \item $u_i \le v_i$ for all $i$,
      \item if $v_i>v_{i+1}$ then $u_i>u_{i+1}$.
    \end{itemize}
  \end{definition}

  \begin{proposition}{\cite[Prop.~5.5]{NadeauTewari}}
    \label{prop:NT}
    The poset $\mathcal P_n$ is a lattice, called (here) the NT lattice.
  \end{proposition}
  In particular, the join of two sequences $u$ and $v$ is the (componentwise) smallest sequence~$w$ whose descent set is contained in the intersection of the descents sets of $u$ and~$v$. For instance, if $u=(4,4,2,2)$ and $v=(4,4,3,1)$ then $w=(4,4,3,3)$. 

  \medskip
  We can  encode an $m$-Dyck path of $\cD_{m,n}$ by a nonincreasing sequence $(u_1, \ldots, u_n)$ of integers, where $u_i$ is the number of $D$ steps that occur after the $i$th large up step $\bU$. For instance, the encoding of the leftmost path of Figure~\ref{figpath} is $(u_1,u_2,u_3)=(6,5,3)$. The encoding sends bijectively $\cD_{m,n}$ on the set of nonincreasing sequences $u$ of length $n$ such that for all $i$,
  \beq\label{condm}
    m(n-i+1) \le u_i \le mn.
  \eeq
  
  Analogously, we encode  a mirrored  $m$-Dyck path of $\cD'_{m,n}$ by the nonincreasing sequence $(u_1, \ldots, u_{mn})$ such that $u_i$ is the number of $\bD$ steps that occur after  the $i$th up step $U$. For instance, the encoding of the rightmost path of Figure~\ref{figpath} is $(u_1, \ldots, u_6)=(3,3,3,2,1,1)$. This encoding sends bijectively $\cD'_{m,n}$ on the set of nonincreasing sequences $u$ of length $mn$ such that for all $i$,
  \beq\label{condmp}
    n-\frac{i-1}m \le u_i \le n.
  \eeq

  \begin{proposition}\label{prop:lattice}
    Let $m,n\geq 1$.    The poset $\LD_{m,n}$ is a lattice, isomorphic to the order induced by the  Nadeau-Tewari order on sequences of  $\mathcal P_n$   satisfying~\eqref{condm}.

   The  poset $\mathbb{D}'_{m,n}$ is a join-semilattice, isomorphic to the order induced by the Nadeau-Tewari order on sequences of $\mathcal P_{mn}$  satisfying~\eqref{condmp}.
 \end{proposition}
 
 \begin{proof}
   The descriptions of $\LD_{m,n}$ and $\LD'_{m,n}$ as induced posets of $\mathcal P_n$ and $\mathcal P_{mn}$ come directly from the characterization of the ascent order given in Proposition~\ref{prop:charac}, and the definitions of   $\LD_{m,n}$ and~$\LD'_{m,n}$ as subposets of $\LD_{mn}$. Let us now address the lattice properties.
   
   In $\mathcal P_n$, any sequence $u$ satisfying~\eqref{condm} satisfies, for the NT order: 
   \[
     u_{\min} :=(mn, m(n-1), \ldots, m) \lqq u \lqq u_{\max}:=(mn, \ldots, mn).
   \]
   This is because $u_{\min}$ has a descent at each place, while $u_{\max}$ has no descent at all. Conversely, any sequence $u$ such that $ u_{\min}  \lqq u \lqq u_{\max}$ satisfies~\eqref{condm}. This means that our encoding sends $\LD_{m,n}$ on the interval $[u_{\min}, u_{\max}]$ of $\mathcal P_n$, and thus  on a lattice.

   Analogously,  any sequence $u$ of $\mathcal P_{mn}$ satisfying~\eqref{condmp} satisfies, for the NT order: 
   \[
    u \lqq u_{\max}:=(n, \ldots, n),
  \]
  because there is no descent in $u_{\max}$. Moreover,  any $v$ such that $u \lqq v \lqq u_{\max}$ also satisfies~\eqref{condmp}.
  Hence our encoding sends $\LD'_{m,n}$ on a union of intervals with maximal element $u_{\max}$, and thus on a join semilattice.
 \end{proof}

\noindent {\bf Examples.} Let us first take $m=1$ and $n=4$, and construct the join in $\LD_n$ of the paths $P=UUDDUUDD$ and $Q= UUDUDDUD$, encoded by $u=(4,4,2,2)$ and $v=(4,4,3,1)$. The join of $u$ and $v$ has already be seen to be $w=(4,4,3,3)$, and thus the join of $P$ and $Q$ is $UUDUUDDD$.

Next let us take $m=2$ and $n=3$, and construct the join in $\LD'_{m,n}$ of $P=UU\bD UUU\bD U\bD$ and $Q= UUUU\bD \bD UU \bD$. The encodings of $P$ and $Q$ are, respectively, $u=(3,3,2,2,2,1)$ and $v=(3,3,3,3,1,1)$, whose join in $\mathcal P_n$ is $w=(3,3,3,3,3,3)$. Hence the join of $P$ and $Q$ in $\LD'_{m,n}$ is~$U^6 \bD^3$. \qee

\section{Intervals in the Nadeau-Tewari poset and classes of the sylvester congruence}
\label{sec:sylvester}

In this section, we exhibit a bijection between (some) intervals of the Nadeau-Tewari lattices~$\cP_n$ and (some) words on the alphabet $\ns^*=\{1,2, \ldots\}$ that \emm avoid the patterns, $aba$ and $acb$ (precise definitions will be given below). These words are known to encode classes of the \emm sylvester congruence, on words, introduced by Hivert et al.~\cite{HNT}. We then specialize our bijection to intervals of $\LD_{m,n}$ and $\LD'_{m,n}$. In particular, intervals of  $\LD'_{m,n}$ are sent bijectively to sylvester classes of $m$-parking functions, considered in~\cite{ThibonNovelli}. 

Let us begin with some definitions. Given a word $w=w_1 \cdots w_n$ on the alphabet $\ns^*$, we denote its length $n$ by $|w|$. We denote by $\{w\}$ the set of letters occurring in $w$, and by $\{\{w\}\}$ the multiset of letters of $w$. We define $\NInc(w)$ (resp. $\NDec(w)$) as the word obtained by reordering the letters of~$w$ in nonincreasing (resp. nondecreasing) order.
For two words~$w$ and $w'$ of the same length, we say that $w\le w'$ if $w_i \le w'_i$ for all $i$. We define $\Low(w)$ as the only nonincreasing word of length $n$ that has the same left-to-right minima as $w$, at the same positions.  In other words, $\Low(w)$ is the largest nonincreasing word (for the above componentwise order) that is smaller than or equal to~$w$. For instance, if $u=6~8~7~4~5~2~3~1~9$ then  ${\Low}(u)=6~6~6~4~4~2~2~1~1$. 

The \emm sylvester congruence, on words is generated by the commutation relations
\[
  ac\cdots b\equiv ca\cdots b, \quad a\leq b<c.
\]
It is known that the words $w=w_1 \cdots w_n$ \emm avoiding the patterns $aba$ and $acb$, form a set of \emm representatives, of sylvester classes~\cite[Sec.~2.7]{ThibonNovelli}: every class contains such a word, and  two distinct words of this type are never congruent.  By \emm pattern avoidance,, we mean that we cannot find $1\le i<j<k\le n$ such that either $w_i=w_k<w_j$ or $w_i<w_k<w_j$.
This congruence arose in connection with binary search trees~\cite{HNT}, and the name \emm sylvester, refers to the forest, \emm silva, in Latin, rather than to the mathematician James Joseph Sylvester.

\begin{lemma}\label{lem:avoid}
  Let $w$ avoid the patterns $aba$ and $acb$. Then $w$ can be uniquely reconstructed from the words $w^{(1)}=\NInc(w)$ and $w^{(2)}=\Low(w)$. 
\end{lemma}
\begin{proof}
By construction, the left-to-right minima of $w$ have the same positions and values as those of $w^{(2)}$. Hence, let us write
  \[
    w= w^{(2)}_{i_1}  z_1 w^{(2)}_{i_2} z_2 \cdots  w^{(2)}_{i_k} z_k
    =w_{i_1}  z_1 w_{i_2} z_2 \cdots  w_{i_k} z_k,
  \]
  where $1=i_1 <i_2 \cdots <i_k$ are the indices of the left-to-right minima of $w$,
  and $|z_j|= i_{j+1}-i_j-1$ for all $j$.
  Starting from $w^{(2)}$ and $w^{(1)}$, we thus need to decide how to arrange the letters of $\{\{w^{(1)}\}\} \setminus \{w_{i_1}^{(2)}, \ldots,w_{i_j}^{(2)} \}$ in the words $z_j$. Obviously, all letters of $z_j$ must be greater than or equal to $w_{i_j}^{(2)}$.

We claim that each  $z_j$ is nondecreasing: indeed, any descent of $z_j$ would give rise to a pattern $aba$ or $acb$, with first letter $a=w_{i_j}^{(2)}$ at position $i_j$, and the other two in $z_j$. So it suffices to determine, inductively on $j$, which letters of $\{\{w^{(1)}\}\}\setminus\{\{w_{i_1}^{(2)} z_1 \cdots z_{j-1} w_{i_j}^{(2)}\}\}$ go into~$z_j$. We claim that they are the $|z_j|$ smallest, among those that are larger than or equal to $w^{(2)}_{i_j}$. Indeed, imagine that we leave out one of them, say $m$. Then the final letter $M$ of $z_j$ is larger than $m$, and $m$ occurs later in $w$. But then the subword $w_{i_j}^{(2)} Mm$ of $w$ has shape  $aba$ or $acb$.
\end{proof}

\begin{example}\label{ex:construction}
  Take $w= 3\ 2\ 2 \ 2 \ 2 \ 5\  1\    1\  1\ 5$, which avoids $aba$ and $acb$. Then $w^{(1)}:=\NInc(w)= 5\ 5\ 3\ 2\ 2\ 2\ 2\ 1\ 1\ 1$ and $w^{(2)}:=\Low(w)=3\ 2\ 2\ 2\ 2\ 2\ 1\ 1\ 1\ 1$. To reconstruct $w$ from these two words, we first collect the left-to-right minima of $w^{(2)}$ and leave them in place:
  \[
    w= 3 \ 2\ \_ \ \_ \ \_ \ \_ 1\_\  \_\  \_.
  \]
The multiset of letters of $w^{(1)}$ that need to fill the gaps is $\{\{5,5,2,2,2,1,1\}\}$.  We fill the first gap~$z_1$, in nondecreasing order, with the $4$ smallest of these  letters that equal at least $2$, namely $2,2,2$ and $5$:
    \[
    w= 3\ 2\ 2 \ 2 \ 2 \ 5 \ 1\_\  \_\  \_.
  \]
  We finally fill the second gap $z_2$ with the remaining letters, $1,1$ and $5$:
   \[
    w= 3\ 2\ 2 \ 2 \ 2 \ 5\  1\    1\  1\ 5.
  \]
  We have recovered $w$.
\end{example}

\begin{figure}[htb]
  \centering
\includegraphics[width=40mm]{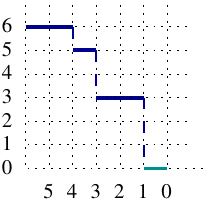}
  \caption{The positive sequence $u=(6,6,5,3,3)\in \cP_5$ is encoded by the word $w=4 \ 3\ 3\ 1\ 1\ 1 \ \in \cW_5$, of length $u_1=6$.}  
  \label{fig:encoding}      
\end{figure}

\begin{remark}
This construction is reminiscent of, but distinct from, the bijection between permutations avoiding the pattern $abc$ and those avoiding $acb$ found in~\cite[Prop.~19]{sim}.
\end{remark}
\medskip
We now consider \emm positive, nonincreasing sequences $u=(u_1, \ldots, u_n)$ of the Nadeau-Tewari poset~$\cP_n$ (Definition~\ref{def:NT}), with $n$ fixed.
These sequences can be encoded bijectively by nonincreasing words $w$
on the alphabet $\{1, \ldots, n\}$ that contain at least one occurrence of the letter~$1$ (Figure~\ref{fig:encoding}): the word $w=W(u)$ has $u_i-u_{i+1}$ occurrences of the letter $n+1-i$, for $1\le i \le n$, where we take $u_{n+1}=0$. Graphically, if we represent $u$ by a path with East steps at heights $u_1, \ldots, u_n,u_{n+1}=0$, joined by South steps, then $w_i$ is the number of East steps after the $i$th South step. It can also be seen as the abscissa of the $i$th South step, if abscissas are numbered from right to left (Figure~\ref{fig:encoding}). Note that $W(u)$ has length $u_1$. We  call it the \emm vertical encoding, of~$u$. Note also that for two sequences $u$ and $v$ in $\cP_n$ such that $u_1=v_1$, we have $u_i\le v_i$ for all $i$ if and only if $W(v)\le W(u)$ componentwise (so the order is  reversed). Of course the sequence~$u$ and the word~$w$ play essentially symmetric roles. We choose to represent $u$ by a sequence and $w$ by a word to keep them distinct. 
We denote by~$\cW_n$ the set of  words on the alphabet $\{1, \ldots, n\}$ that contain at least one occurrence of the letter~$1$.  Nonincreasing words of $\cW_n$ are in bijection with positive sequences of $\cP_n$.

\begin{proposition}
 Fix $n\ge 1$ and  let $w \in \cW_n$. Then  $\NInc(w)$ and $\Low(w)$ are nonincreasing words of the same length in $\cW_n$. Let $u$ (resp. $v$) be the positive sequence of $\cP_n$ such that $W(u)=\NInc(w)$ (resp. $W(v)=\Low(w)$). Then $u_1=v_1=|w|$ and $u\lqq v$ in the Nadeau-Tewari poset~$\cP_n$. Define $\Phi_n(w)$ to be the interval $[u,v]$.

  The restriction of $\Phi_n$ to words of $\cW_n$ avoiding the patterns $aba$ and $acb$ is a bijection between these words and  intervals $[u,v]$ of positive words  such that $u_1=v_1$ in the Nadeau-Tewari poset~$\cP_n$. 
\end{proposition}
\begin{proof}
  Let us first observe that, for two positive sequences $u$ and $v$ in $\cP_n$ such that $u_1=v_1$,  the condition $u\lqq v$ translates in terms of the words $w^{(1)}:=W(u)$ and $w^{(2)}:=W(v)$ as follows: $w^{(2)}\le w^{(1)}$ (componentwise), and the set $\{w^{(2)}\}$ of letters of $w^{(2)}$ is included in $\{w^{(1)}\}$. 

  Now let $w=w_1 \cdots w_M$ be a word of length $M$ in $\cW_n$, and take $w^{(1)}=\NInc(w)$ and $w^{(2)}=\Low(w)$. Write $w^{(k)}= w^{(k)}_1\cdots w^{(k)}_M$ for $k=1, 2$. Given that the letters of $w^{(2)}$ are the values of the left-to-right minima of $w$, while $w^{(1)}$ is just a reordering of $w$, we have $\{w^{(2)}\}\subset \{w^{(1)}\}$. Let us now prove that $w^{(2)}\leq w^{(1)}$ componentwise (note that for the moment we do not assume that $w$ avoids any pattern).
  For any letter $a \in \ns$, we have, by definition of $w^{(1)}=\NInc(w)$:
  \[
    w^{(1)}_i\ge a \Leftrightarrow \,\sharp \{j: w_j \ge a\} \ge i.
  \]
  Hence $w^{(2)}_i\le w^{(1)}_i$ if and only if at least $i$ letters of $w$ are  larger than or equal to~$w^{(2)}_i$. Assume that $m\le i <m'$, where $m$ and $m'$ are the positions of two consecutive left-to-right minima of~$w$. Then $w^{(2)}_i=w_m$, and by definition of $w^{(2)}=\Low(w)$, the letters $w_1, \ldots, w_{m'-1}$ are larger than or equal to~$w_m$. Since there are $m'-1$ of them, and $i\le m'-1$, this proves that $w^{(2)}_i\le w^{(1)}_i$. We conclude that the positive sequences $u$ and $v$ of $\cP_n$ given by $W(u)=w^{(1)}$ and $W(v)=w^{(2)}$ form an interval for the NT order. They also satisfy $u_1=v_1=M$.
  These properties hold in particular when $w$ avoids the two forbidden patterns. Moreover, the map $\Phi_n$ is injective on those words, by Lemma~\ref{lem:avoid}.

Let us now prove surjectivity. Let $u\lqq v$ in $\cP_n$, with $u$ and $v$ positive and ${u_1=v_1:=M}$. Let $w^{(1)}=W(u)$ and $w^{(2)}=W(v)$. As already observed,  $u\lqq v$  means that $w^{(2)}\le w^{(1)}$ (componentwise) and $\{w^{(2)}\}\subset \{w^{(1)}\}$. The proof of Lemma~\ref{lem:avoid} tells us how to reconstruct a word $w$ that avoids the two patterns and satisfies $\Phi(w)=(u,v)$ --- if such a word exists! So let us try to apply this construction to  $w^{(1)}$ and $w^{(2)}$: we keep the left-to-right minima of $w^{(2)}$, denoted~$w_{i_j}^{(2)}$, in place, and fill the gaps $z_j$ with the remaining letters of $w^{(1)}$. Let us first explain that the construction succeeds, that is, produces a word $w$.  We will then explain why $w$ does not contain the forbidden patterns. The only thing that could go wrong in the construction of $w$ is that we could run out of letters when filling,  with letters of $\{\{w^{(1)}\}\} \setminus \{\{w_{i_1}^{(2)} z_1 \cdots z_{j-1} w_{i_j}^{(2)}\}\}$, the gap following the $j$th  left-to-right minimum of $w$. This only happens if the number of letters of $w^{(1)}$ that are larger than or equal to $w_{i_j}^{(2)}$ is less than $i_{j+1}-1$ (with $i_{j+1}=M+1$ if $w_{i_j}^{(2)}$ is the rightmost left-to-right minimum). However, the number of such letters in $w^{(2)}$ is precisely $i_{j+1}-1$, and $w^{(2)}\le w^{(1)}$ componentwise: hence there are at least $i_{j+1}-1$ letters larger  than or equal to $w_{i_j}^{(2)}$ in  $w^{(1)}$, and one is never stuck in the construction of $w$. Finally, $w$ avoids $aba$ and $acb$: indeed, if there was one of these patterns in $w$, there would be one where the first $a$ is occurs at a left-to-right minimum, hence  $a=w_{i_j}^{(2)}$: but the construction has been designed to guarantee that there is no $w_k>w_\ell\ge a$ with $i_j<k<\ell$.  Hence the map~$\Phi_n$ is surjective, and the proposition is proved.  
\end{proof}

\begin{figure}[htb]
  \centering
\includegraphics[width=40mm]{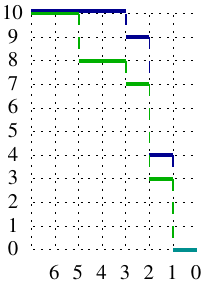}
  \caption{The sequences $u=(10,10,8,8,7,3)$ and $v= (10,10,10,10,9,4)$ form an interval in $\cP_6$ such that $u_1=v_1$. Their vertical encodings are $w^{(1)} = 5\ 5\ 3\ 2\ 2\ 2\ 2\ 1\ 1\ 1$ and $w^{(2)}=3\ 2\ 2\ 2\ 2\ 2\ 1\ 1\ 1\ 1$,  respectively, and  the associated pattern avoiding word is  $w= 3\ 2\ 2 \ 2 \ 2 \ 5\  1\    1\  1\ 5$ (Example~\ref{ex:bij}).}  
  \label{fig:ex-bij}       
\end{figure}

\begin{example}\label{ex:bij}
  Take  $n=6$, and choose $w= 3\ 2\ 2 \ 2 \ 2 \ 5\  1\    1\  1\ 5$ as in Example~\ref{ex:construction}. Then  $w^{(1)}:=\NInc(w)= 5\ 5\ 3\ 2\ 2\ 2\ 2\ 1\ 1\ 1$ and $w^{(2)}:=\Low(w)=3\ 2\ 2\ 2\ 2\ 2\ 1\ 1\ 1\ 1$.  The positive nonincreasing sequences $u$ and $v$ of $\cP_n$ such that $W(u)=w^{(1)}$ and $W(v)=w^{(2)}$ are $u=(10,10,8,8,7,3)$ and $v= (10,10,10,10,9,4)$, and they form an interval in $\cP_n$ (Figure~\ref{fig:ex-bij}).
\end{example}

It follows from the above proposition that, given a positive sequence $u^{(0)}$ of  $\cP_n$, the map $\Psi_n:=\Phi_n^{-1}$ sends  bijectively  intervals $[u,v]$ of $\cP_n$ such that $u_1=v_1=u^{(0)}_1$
and $u_i\ge u_i^{(0)}$ for all~$i\le n$ onto
\begin{itemize}
\item  words $w$ of $\cW_n$ avoiding $aba$ and $acb$ and satisfying  $\NInc(w)\le W(u^{(0)})$, componentwise,
  \item or equivalently sylvester classes of words $w$ of $\cW_n$ such that  $\NInc(w)\le W(u^{(0)})$.
\end{itemize}
Returning to Proposition~\ref{prop:lattice}, we can use this to exhibit families of words (or sylvester classes) in bijection with intervals of $\LD_{m,n}$ and $\LD'_{m,n}$.

\begin{corollary}\label{cor:int-sylvester} Upon encoding $m$-Dyck paths and mirrored $m$-Dyck paths by nonincreasing sequences (see Section~\ref{sec:lattice}), the maps $\Psi_n$ and $\Psi_{mn}$ induce respectively  bijections between
  \begin{itemize}
  \item intervals of $\LD_{m,n}$  and sylvester classes of words $w$ of $\cW_n$ of length $mn$
    such that $\NInc(w)\le n^m (n-1)^m \cdots 1^m$, where $i^m$ represents  a sequence of $m$ copies of $i$;
  \item intervals of $\LD'_{m,n}$  and sylvester classes of words $w$ of $\cW_{mn}$ of length $n$    such that $\NInc(w)\le ((n-1)m+1)\cdots (2m+1)(m+1) 1$.
  \end{itemize}
\end{corollary}
Note that the two statements coincide when $m=1$.
The second one explains why we found the sequences $\left(g'_m(n)\right)_{n>0}$ counting intervals of $\LD'_{m,n}$, for $1\le m \le 5$, in the OEIS: they appear at the very end of~\cite{ThibonNovelli} in Table~21, as counting sylvester classes of $m$-parking functions (of size $n$). These functions are defined in~\cite[Sec.~6.1]{ThibonNovelli} as positive words $w$ of length $n$ such that $\NDec(w)\le 1 (m+1) \cdots ((n-1)m+1)$, which is equivalent to saying that $\NInc(w)\le((n-1)m+1)\cdots (m+1) 1$ as above. The sylvester classes involved in the first part of the corollary do not seem to have been considered so far.

\begin{example} Let us take $m=1$ and $n=4$, and the interval $[P,Q]$ given by $P=UDUUDDUD$ and $Q=UUUDDUDD$. We encode $P$ and $Q$ by nonincreasing sequences $u$ and $v$ as described in Section~\ref{sec:lattice}, with $u=(4,3,3,1)$ and $v=(4,4,4,2)$. Note that $u\lqq v$ for the NT order, as used in the proof of Proposition~\ref{prop:lattice}. The vertical encodings of $u$ and~$v$ are $W(u)=4\ 2\ 2\ 1$ and $W(v)=2\ 2\ 1\ 1$. The unique word $w$ avoiding $aba$ and $acb$ such that $\NInc(w)=W(u)$ and $\Low(w)=W(v)$ is $w=2\ 2\ 1\ 4$ and is the canonical representative of the sylvester class associated with $[P,Q]$. 
\end{example}

\begin{remark}
    Words avoiding the patterns $aba$ and $acb$ have recently been considered for their own sake. See~\cite[Sec.~4.7]{testart}, where a formula is given for the number of such words of length $n$ on a $k$ letter alphabet. However, the underlying recursive construction does not seem compatible with the conditions of Corollary~\ref{cor:int-sylvester}.
\end{remark}

\section{Recursive construction of intervals in $\mathbb{D}_{m,n}$ and $\mathbb{D}'_{m,n}$}
\label{sec:recursive}

In this section, we describe recursive constructions of the  intervals of the posets $\mathbb{D}_{m,n}$ and $\mathbb{D}'_{m,n}$, and convert them into functional equations for the associated \gfs.

\subsection{Intervals in $\mathbb{D}_{m,n}$ }

\begin{lemma}\label{lem:final-insert}
  Let $m, n\ge 1$, and let $[P',Q']$ be an interval of $\LD_{m,n}$. Delete in $P'$ and $Q'$ the final (large) peak $\bU D^m$, to obtain two paths $P$ and $Q$ of $\LD_{m,n-1}$. Then they form an interval, that is, $P\le Q$.

  Conversely, start from an interval $[P,Q]$ in $\LD_{m,n-1}$. Let $a\le b$ be the lengths of the final descents of~$P$ and $Q$, respectively. Inserting  final large peaks $\bU D^m$ in the final descents of $P$ and $Q$, starting at respective heights $a' \in \llbracket 0,a\rrbracket$ and $b' \in \llbracket 0,b\rrbracket$, yields an interval $[P',Q']$ if and only if $a'\le b'$ and $b'$ is maximal as soon as $a'$ is maximal. That is, $a'=a$ implies $b'=b$.
\end{lemma}
\begin{proof}
  The first part directly follows from Proposition~\ref{prop:charac}: deleting the final large peak preserves the relation ``lying below'', and subtracts $m$ from the final part of the ascent composition (parts of size $0$ are then discarded), thus preserving refinement. 

  Conversely, inserting   new large peaks  in the final descents of $P$ and $Q$, at heights $a'$ and~$b'$ respectively, preserves the relation ``lying below'' if and only if $a'\le b'$. This operation adds to~$c(P)$ a new part equal to $m$ if $a'<a$, and  adds $m$ to the final part of $c(P)$ if $a'=a$. An analogous statement holds for $Q$. Since we need $c(P')$ to refine $c(Q')$, these observations imply the final condition of the lemma.
\end{proof}

Let $G_m(t;x,y)$  be the \gf\ of intervals $[P,Q]$ in the ascent lattices $\mathbb{D}_{m,n}$, counted by the size $n>0$ (variable $t$), the length of the final descent of $P$ (variable $x$) and of $Q$ (variable $y$). Then $G_m(t;x,y)$  are formal power series in $t$ whose coefficients are polynomials in $x$ and $y$. For $m=1,2$, they start as follows:
\[
  G_1(t;x,y)
  =x yt
  + x y \left(x y +y +1\right) t^{2}
  +x y \left(x^{2} y^{2}+2 x \,y^{2}+2 x y +2 y^{2}+3 y +3\right) t^{3}  + \LandauO(t^4),
\]
\begin{multline*}
  G_2(t;x,y)= x^2 y^2t  + x^2 y^2  \left(x^2 y^2+x y^2+xy +y^2+y +1\right) t^{2}  +
  x^2 y^2 \left(x^4y^4+2x^3y^4+2x^3y^3+3x^2y^4\right.\\
  \left. +4x^2y^3+3xy^4+4x^2y^2+5xy^3+3y^4+6xy^2+5y^3+6xy+6y^2+6y+6\right)t^{3}  + \LandauO(t^4).
  \end{multline*}
Clearly, the exponent of $y$ is at least equal to the exponent of $x$ in each monomial of $G_m(t;x,y)$, and in addition all terms in $G_m$ are multiples of $tx^m y^m$. So we define a new series  $Q_m(t;x,y)$ by
\beq\label{HQ}
  G_m(t;x,y)=t x^m y^m Q_m(t;xy,y),\quad \mbox{ or equivalently } \quad Q_m(t;x,y)=G_m(t;x/y,y)/(tx^m).
\eeq

\begin{proposition}\label{prop:eq-up}
 For $m\ge 1$,  the \gf\ $G_m(t;x,y)\equiv G_m(x,y)$ is characterized by the following functional equation:
  \[
    G_m \! \left(x , y\right)=t x^m y^m +tx^m y^{m}  G_m \! \left(x , y\right)
    +t x^m \,y^{m+1}\, \frac{ G_m \! \left(x , y\right)-G_m \! \left(1, y\right)}{\left(x -1\right) \left(y -1\right)}
    -tx^m y^m \, \frac{ G_m \! \left(x y , 1\right)-G_m \! \left(1, 1\right)}{\left(y -1\right) \left(x y -1\right)}
    .\]
  Equivalently, $Q_m(t;x,y)\equiv Q_m(x,y)$ is characterized by
  \[
    Q_m \! \left(x , y\right)=1  +tx^m Q_m \! \left(x , y\right)
    +t  y^2\, \frac{x^m Q_m \! \left(x , y\right)-y^mQ_m \! \left(y, y\right)}{\left(x -y\right) \left(y -1\right)}
    -t  \,\frac{x^m Q_m \! \left(x  , 1\right)-Q_m \! \left(1, 1\right)}{\left(x -1\right) \left(y  -1\right)}.
  \]
\end{proposition}

\begin{proof}
  We construct  intervals of $\LD_{m,n}$  recursively on the size $n$, starting from the only interval of size $1$, namely $[\bU D^m , \bU D^m]$, and inserting a final peak as in Lemma~\ref{lem:final-insert}. We use the notation of this lemma. The rule that describes the final descent lengths of paths $P'\le Q'$ obtained from  $[P,Q]\in \LD_{m,n-1}$ in terms of the final descent lengths $a$ and $b$ of $P$ and $Q$ is
  \beq
  (a,b)  \rightarrow \left\{
    \begin{array}{ll}
    (m+a', m+b'), & \text{ for } \ 0\le a' <a \ \text{ and } \ a'\le b' \le b,  \\
      (m+a, m+b).
       \end{array}
    \right.
    \label{rr}
  \eeq
  In other words, intervals of  $\LD_{m,n}$ are in bijection with nodes at height $(n-1)$ in the \emm generating tree, having root $(m,m)$ and the above rewriting rule (see for instance~\cite{west-GT,bousquet-motifs} for the hopefully intuitive notion of generating trees). In terms of \gfs, if we write $G_m(x,y)=\sum_{0\le a\le b} G_{m,a,b} x^a y^b$, so that $G_{m,a,b} $ is the series (in $t$) counting intervals with final descent lengths~$a$ and $b$, the above construction gives:
  \allowdisplaybreaks
  \begin{align*}
    G_m(x,y)&=tx^m y^m + t\sum_{0\le a\le b} G_{m,a,b}\left( x^{m+a} y^{m+b}
              +\sum_{a'=0}^{a-1} \sum_{b'=a'}^{b} x^{m+a'} y^{m+b'}\right)\\
            & = tx^m y^m + tx^m y^m G_m(x,y)+ tx^m y^m\sum_{0\le a\le b} G_{m,a,b} \sum_{a'=0}^{a-1}x^{a'}\cdot \frac{y^{b+1}-y^{a'}}{y-1}\\
            & = tx^m y^m + tx^m y^m G_m(x,y)+ tx^m y^m\sum_{0\le a\le b} G_{m,a,b} \frac{x^a-1}{x-1} \frac{y^{b+1}}{y-1}\\
    & \hskip 60mm -tx^m y^m\sum_{0\le a\le b} G_{m,a,b} \frac{x^ay^a-1}{(xy-1)(y-1)} ,    
  \end{align*}
  which is the announced equation for $G_m$. It is straightforward to convert it into an equation for~$Q_m$, using~\eqref{HQ}.
\end{proof}

The rewriting rule~\eqref{rr}
  also describes the construction of certain lattice walks confined to the first quadrant of the plane, which are thus in bijection with ascent intervals. This will be combined in Section~\ref{sec:asympt} to general probabilistic results on quadrant walks to obtain the asymptotic result of Proposition~\ref{prop:Dmn}.

\begin{corollary}\label{cor:bij-mn}
  Let $m,n\ge 1$. There is a bijection between  intervals in $\LD_{m,n}$ and walks  in the quarter plane $\ns^2$ that start from $(0,0)$ and consist of $n-1$ steps taken from the following subset of $\zs^2$:
  \[
    \cS_m= \{(m,0)\} \cup \left\{ (\delta_x, \delta_y): \delta_x<m \text{ and } \delta_x+ \delta_y \le m \right\}.
  \]
  More precisely, a walk ending at $(i,j)$ corresponds to an interval $[P,Q]$ where $P$ and $Q$ have final descent lengths $m+i$ and $m+i+j$, respectively.

  This also yields a bijection between intervals in $\LD_{m,n}$ and quadrant walks of length $n$ starting and ending at $(0,0)$, still taking their steps in $\cS_m$.  We call such walks \emm excursions,.
\end{corollary}

\begin{proof}
  There are two (essentially equivalent) ways of proving this result.

  The first one starts from the rewriting rule~\eqref{rr} and rewrites the label $(a,b)$ as $(m+i, m+i+j)$, so that the conditions $m\le a\le b$ become $i,j \ge 0$. The transformed tree has root $(0,0)$ and rewriting rule
  \beq
  (i,j) \rightarrow
  \left\{
    \begin{array}{ll}
     (k, \ell), &  \text{ for }\ 0\le k < m+i \ \text{ and }\  0\le \ell \le m+i+j-k, \\
     (m+i,j).    
    \end{array}
     \right. \label{ij}
  \eeq
  In other words, from the point $(i,j)$ in the nonnegative quadrant $\ns^2$, we can move to another point $(k, \ell)$ of the quadrant by appending a step $(\delta_x, \delta_y)=(k-i, \ell-j)$ satisfying
  \beq
    (\delta_x, \delta_y)= (m,0) \quad \text{or} \quad \left( \delta_x<m \text{ and } \delta_x+ \delta_y \le m\right), \label{steps}
  \eeq
  which is indeed the collection of steps $\cS_m$. The second statement follows by adding to a walk of length $n-1$ ending at $(i,j)$ the step $(-i,-j)$, which is indeed in $\cS_m$.

  The second way to prove the corollary is to start from the description of the quadrant walk and to recover the equation on $Q_m$   obtained in Proposition~\ref{prop:eq-up}. Since handling quadrant walks with infinitely many steps is uncommon in the quadrant walks literature, let us do this. Let us denote by $Q(t;x,y)\equiv Q(x,y)=\sum_{i,j\ge 0} Q_{i,j} x^i y^j$ the \gf\ of quadrant walks with steps in $\cS_m$ (starting from $(0,0)$), counted by the length (i.e., the number of steps; variable $t$), and the coordinates of the endpoint (variables $x$ and $y$). We construct walks step-by-step, using the description~\eqref{ij} of the endpoint  rather than the description~\eqref{steps} of the steps:
  \begin{align*}
    Q(x,y)&=1+ t x^m Q(x,y) + t \sum_{i, j \ge 0} Q_{i,j} 
            \left(\sum_{k=0}^{m+i-1}  \sum_{\ell=0}^{m+i+j-k} x^k y^\ell\right)\\
    &=1+ t x^m Q(x,y) + t \sum_{i, j \ge 0} Q_{i,j} 
      \left(\sum_{k=0}^{m+i-1}   x^k \cdot \frac{ y^{m+i+j-k+1}- 1}{y-1}\right)\\
&=1+ t x^m Q(x,y) + \frac t{y-1} \sum_{i, j \ge 0} Q_{i,j} 
                                                                                    \left( y^{m+i+j+1} \frac{(x/y)^{m+i} -1}{x/y-1}
                                                                                    - \frac{x^{m+i}-1}{x-1}\right),
  \end{align*}
  which gives the equation obtained in Proposition~\ref{prop:eq-up} for the series $Q_m$.
\end{proof}

\subsection{Intervals in $\mathbb{D}'_{m,n}$ }
\label{sec:recp}
Let us now turn our attention to mirrored $m$-Dyck paths. First, observe that the final descent $\bD^h$ of such a path is not necessarily preceded by at least $m$ up steps, so there is not always a final (large) peak. It makes sense to consider instead the \emm first, (large) peak $U^m \bD$, which always exists. We will use a new parameter on intervals: when $P\le Q$, with associated  ascent compositions $c(P)=(c_1,c_2, \ldots)$ and $c(Q)=(d_1, d_2, \ldots)$, we know that $c(P)$ refines $c(Q)$. We define  $r(P,Q)$ to be the unique integer $r$ such that $c_1 + \cdots + c_r=d_1$. 

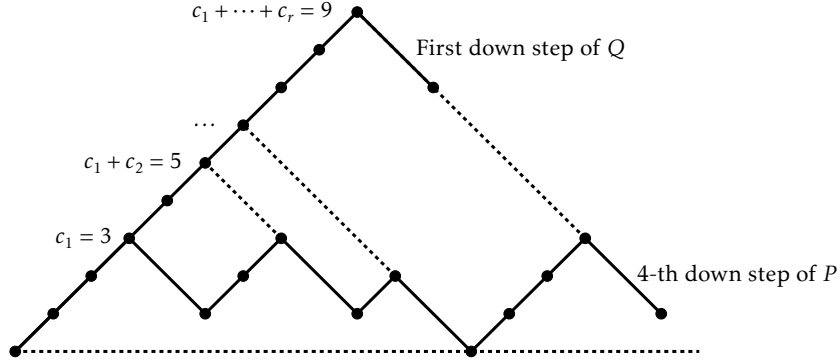
\begin{figure}[htb]
   \centering\scalebox{0.5}{
\begin{tikzpicture}
\draw[line width = 0.8mm, dashed] (0,0) -- (18,0);
\draw[line width = 0.8mm] (0,0) -- (1,1)--(2,2)--(3,3)--(4,4)--(5,5)--(6,6)--(7,7)--(8,8)--(9,9)--(11,7);
\filldraw[black] (0,0) circle (4.pt);
\filldraw[black] (1,1) circle (4.pt);
\filldraw[black] (2,2) circle (4.pt);
\filldraw[black] (3,3) circle (4.pt);
\filldraw[black] (4,4) circle (4.pt);
\filldraw[black] (5,5) circle (4.pt);
\filldraw[black] (6,6) circle (4.pt);
\filldraw[black] (7,7) circle (4.pt);
\filldraw[black] (8,8) circle (4.pt);
\filldraw[black] (9,9) circle (4.pt);
\filldraw[black] (11,7) circle (4.pt);
\draw (13.3,8) node {\LARGE{First down step of $Q$}}; 
\draw[line width = 0.8mm] (0,0) -- (1,1)--(2,2)--(3,3)--(5,1)--(6,2)--(7,3)--(9,1)--(10,2)--(12,0)--(13,1)--(14,2)--(15,3)--(17,1);
\filldraw[black] (5,1) circle (4.pt);
\filldraw[black] (6,2) circle (4.pt);
\filldraw[black] (7,3) circle (4.pt);
\filldraw[black] (9,1) circle (4.pt);
\filldraw[black] (10,2) circle (4.pt);
\filldraw[black] (12,0) circle (4.pt);
\filldraw[black] (13,1) circle (4.pt);
\filldraw[black] (14,2) circle (4.pt);
\filldraw[black] (15,3) circle (4.pt);
\filldraw[black] (17,1) circle (4.pt);
\draw[line width = 0.8mm, dashed] (15,3) -- (9,9);
\draw[line width = 0.8mm, dashed] (10,2) -- (6,6);
\draw[line width = 0.8mm, dashed] (7,3) -- (5,5);
\draw (1.8,3) node {\LARGE $c_1=3$};
\draw (3.1,5) node {\LARGE{$c_1+c_2=5$}};
\draw (5,6) node {\LARGE $\ldots$};
\draw (6.5,9) node {\LARGE $c_1+ \cdots + c_r=9$};
\draw (19,2) node {\LARGE{$4$-th down step of $P$}};
\end{tikzpicture}}
\caption{An example where  $r(P,Q)=4$ for an interval $[P,Q]$ in  $\mathbb{D}'_{2,n-1}$. A first peak can be inserted in $Q$ at heights $b'=0, 1, \ldots, c_1, c_1+c_2, \ldots, c_1+\cdots + c_r$. If $b'<c_1$, then the peak inserted in $P$ must start at height $b'$ as well.}  
\label{fig:newstat}
\end{figure}

\begin{lemma}\label{lem:first-insert}
  Let $m, n\ge 1$, and let $[P',Q']$ be an interval of $\LD'_{m,n}$. Delete in $P'$ and $Q'$ the first (large) peak $U^m\bD$, to obtain two paths $P$ and $Q$ of $\LD'_{m,n-1}$. Then they form an interval, that is, $P\le Q$.

  Conversely, start from an interval $[P,Q]$ in $\LD'_{m,n-1}$. Let
   $c(P)=(c_1,c_2, \ldots)$ and $c(Q)=(d_1, d_2, \ldots)$ be the corresponding ascent compositions. The length of the first ascent in $P$ (resp. $Q$) is thus $c_1$ (resp.~$d_1$). Let $r=r(P,Q)$.
  Inserting a first large peak $U^m\bD$ in the first ascents of $P$ and $Q$, starting at heights $a'\in \llbracket 0, c_1\rrbracket $ and $b'\in \llbracket 0, d_1\rrbracket$ respectively, yields an interval if and only if
  \begin{itemize}
  \item either $b'$ is one of $c_1, c_1+c_2, \ldots, c_1+\cdots + c_r=d_1$ and $a'$ is any element of $\llbracket 0, c_1\rrbracket $,
    \item or $0\le a'=b'<c_1$.
  \end{itemize}
\end{lemma}
\begin{proof}
  The first part  follows  from Proposition~\ref{prop:charac}: when deleting the first peaks, one just subtracts~$m$ from the first part in the ascent compositions of $P'$ and $Q'$.

  Conversely, inserting   new large peaks  in the first ascents of $P$ and $Q$ at heights $a'$ and~$b'$ respectively, preserves the relation ``lying below'' if and only if $a'\le b'$. We want to choose $a'$ and~$b'$ so that the resulting paths $P'$ and $Q'$ are such that $c(P')$ refines $c(Q')$. Recall that this means that every descent of $Q'$ lies on the same line (of slope $-1$) as a descent of $P'$. This gives the conditions stated in the lemma; see Figure~\ref{fig:newstat}. 
\end{proof}


Let   $G'_m(t;x,y)$ be the \gf\ of intervals $[P,Q]$ in the ascent lattice $\mathbb{D}'_{m,n}$, counted by the size $n$ (variable $t$), the length of the first ascent of $P$ (variable $x$) and the statistic $r(P,Q)$ (variable $y$).  Then $G'_m(t;x,y)$ is a  formal power series in $t$ whose coefficients are polynomials in $x$ and $y$. For $m=1, 2$, they start as follows:
\begin{multline*}
   G'_1(t;x,y)=
   x y  t +xy \left(x +y +1\right)   t^{2}+x y \left(x^{2}+2 x y +y^{2}+3 x +3 y +3\right) t^{3}\\
   +x y \left(x^{3}+3 x^{2} y +3 x \,y^{2}+y^{3}+6 x^{2}+10 x y +6 y^{2}+13 x +13 y +13\right) t^{4}+\LandauO \left(t^{5}\right),
\end{multline*}
 \begin{align*} 
  G'_2(t;x,y)&= x^2yt
  + x^2 y\left(x^2 +xy+x+y +1\right) t^{2}+\\
&  +x^2 y\left(x^4+2x^3y+x^2y^2+3x^3+4x^2y+2xy^2+5x^2+5xy+2y^2+5x+5y+5\right) t^{3}  + \LandauO(t^4).
 \end{align*}
 Note that $G'_1(t;x,y)$ does not coincide with $G_1(t;x,y)$, except at the point $x=y=1$, since these series record different statistics. We will see below that $G'_1(t;x,y)$ is symmetric in $x$ and~$y$, a result that calls for a direct combinatorial explanation.

 Since all monomials are multiples of $tx^m y$, we introduce the series $Q'_m(t;x,y)$ defined by
 \beq\label{HQp}
   G'_m(t;x,y)=tx^m y Q'_m(t;x,y).
 \eeq
 
\begin{proposition}\label{prop:eq-down}
  The \gf\ $G'_m(t;x,y)\equiv G'_m(x,y)$ is characterized by the following functional equation:
  \[
    G'_m \! \left(x , y\right)=t x^m y+tx^my\frac{G'_m\! \left(x , y\right)-G'_m \! \left(x, 1\right)}{y-1}     +t x^m \,y^{2}\, \frac{ G'_m \! \left(x , y\right)-G'_m \! \left(1, y\right)}{\left(x -1\right) \left(y -1\right)}
    -tx^my \, \frac{ G'_m \! \left(x , 1\right)-G'_m \! \left(1, 1\right)}{\left(x -1\right)\left(y -1\right) }
    .\]
  Equivalently, $Q'_m(t;x,y)\equiv Q'_m(x,y)$ is characterized by
  \[
    Q'_m(x,y)=1+tx^m\frac{yQ'_m\! \left(x , y\right)-Q'_m \! \left(x, 1\right)}{y-1}
    +t y^{2}\, \frac{ x^m Q'_m \! \left(x , y\right)-Q'_m \! \left(1, y\right)}{\left(x -1\right) \left(y -1\right)}
    -t\, \frac{ x^m Q'_m \! \left(x , 1\right)-Q'_m \! \left(1, 1\right)}{\left(x -1\right)\left(y -1\right) }
    .\]
\end{proposition}

\begin{remark}
  When $m=1$, the equation for $G'_1$ can be rewritten as
 \[
   G'_1 \! \left(x , y\right)=t x y+txy\frac{(x+y-1)G'_1\! \left(x , y\right)-xG'_1 \! \left(x, 1\right)-yG'_1 \! \left(1, y\right)+G(1,1)}{\left(x -1\right)\left(y -1\right) }
   .
 \]
 It is symmetric in $x$ and $y$, and characterizes $G'_1(x,y)$: we conclude that this series is symmetric in~$x$ and $y$. Since the equation reflects a (symmetric) recursive construction of intervals, it can be used to define a recursive involution on ascent intervals that exchanges the length of the first ascent of the bottom path and the parameter $r(\cdot, \cdot)$. See Section~\ref{sec:final} for details. Note that another symmetry property was discovered in Tamari intervals~\cite{mbm-fusy-preville}, and explained in~\cite{chapoton-chatel-pons,pons-involution}.
\end{remark}

\begin{proof}[Proof of Proposition~\ref{prop:eq-down}]
  We construct  intervals of $\LD'_{m,n}$  recursively in the size $n$, starting from the only interval of size $1$, namely $[U^m \bD , U^m\bD]$, and inserting  peaks as in Lemma~\ref{lem:first-insert}. We use the notation of this lemma. Let us examine  the value of the parameter $r(P',Q')$, for  paths $P'\le Q'$ obtained from  $[P,Q]\in \LD'_{m,n-1}$:
  \begin{itemize}
  \item if $b'=c_1+ \cdots + c_s$, then $r(P',Q')=s$ if $a'=c_1$, and $r(P',Q')=s+1$ otherwise,
    \item if $b'=a'<c_1$, then $r(P',Q')=1$.
  \end{itemize}
Hence, if we record for each interval $[P,Q]$ the pair $(c_1(P), r(P,Q))$, we obtain a generating tree with root labelled $(m,1)$ and rewriting rule
 \begin{equation} \label{rr-bis} 
(a,r)  \rightarrow \left\{ 
  \begin{array}{ll}
   (m+a, s), & \text{ for } 1\le s \le r, \\
      (m+a', s+1), & \text{ for } 0\le a' < a \text{ and } 1\le s \le r, \\
       (m+a', 1), & \text{ for } 0\le a'<a.   
  \end{array}
       \right. 
  \end{equation}
  Note that the last two lines can be merged into a single one, by allowing  $s$ to be $0$.
  
  In terms of \gfs, if we write $G'_m(x,y)=\sum_{ a,r} G'_{m,a,r} x^a y^r$, the above construction gives:
  \begin{align*}
    G'_m(x,y) & = t x^m y + t \sum_{a,r} G'_{m,a,r}\left( x^{m+a} \sum_{s=1}^r y^s
                + \sum_{a'=0}^{a-1} \sum_{s=0}^r x^{m+a'}  y^{s+1}
              \right)\\
    & = t x^m y + t x^m\sum_{a,r} G'_{m,a,r}\left( x^{a}\cdot \frac{y^{r+1}-y}{y-1} + \frac{x^a-1}{x-1} \cdot \frac{y^{r+2}-y}{y-1}\right),
  \end{align*}
  which gives the first equation of Proposition~\ref{prop:eq-down}. It is straightforward to convert it into an equation for $Q'_m$, using~\eqref{HQp}.
\end{proof}

\begin{corollary} \label{cor:bij-mnp}
  Let $m,n\ge 1$. There is a bijection between  intervals in $\LD'_{m,n}$ and walks  in the quarter plane $\ns^2$ that start from $(0,0)$ and consist of $n-1$ steps taken from
  \[
    \cS'_m= \big(\{m\} \times \llparenthesis -\infty, 0\rrbracket\big)
    \cup \big( \llparenthesis -\infty, m-1\rrbracket \times \llparenthesis -\infty , 1\rrbracket \big),
  \]
  where $\llparenthesis -\infty, a\rrbracket:= \zs \cap (-\infty, a]$.
  Moreover, a walk ending at $(i,j)$ corresponds to an interval $[P,Q]$ where $P$  has a first ascent of length $m+i$ and $r(P,Q)=1+j$.

  This also yields a bijection between intervals in $\LD'_{m,n}$ and quadrant walks of length $n$ starting and ending at $(0,0)$, still taking their steps in $\cS'_m$.  
\end{corollary}
\begin{proof}
Let us start from the rewriting rule~\eqref{rr-bis} and rewrite the label $(a,r)$ as $(m+i, 1+j)$, so that the conditions $m\le a, 1\le r$ become $i,j \ge 0$. The transformed tree has root $(0,0)$ and rewriting rule
  \beq
  (i,j) \rightarrow
  \left\{
    \begin{array}{ll}
     (m+i, \ell) & \text{ for }\  0\le \ell   \le j\\
     (k,\ell ) & \text{ for } \ 0\le k <i+m \ \text{ and }\ 0\le \ell \le 1+j.
    \end{array}
    \right.
  \eeq
  In other words, from the point $(i,j)$ in the quadrant, we can move to another point $(k, \ell)$ of the quadrant by appending a step $(\delta_x, \delta_y)=(k-i, \ell-j)$ satisfying
  \beq
    \left(\delta_x= m \text{ and }  \delta_y\le 0\right) \quad \text{or} \quad \left( \delta_x<m \text{ and }  \delta_y \le 1\right), \label{steps-bis}
  \eeq
  which is indeed the collection of steps $\cS'_m$. The second statement follows by adding to a walk of length $n-1$ ending at $(i,j)$ the step $(-i,-j)$, which is indeed in $\cS'_m$.
\end{proof}

There is an alternative description of the series $G'_m(1,1)$ in terms of quadrant walks  involving \emm finitely many weighted steps,, or equivalently steps from a \emm finite multiset,.

\begin{corollary}\label{cor:finite}
Let $m, n \ge 1$.  The number of  intervals in the poset $\LD'_{m,n}$ is the number of quadrant walks of length $n-1$ starting and ending at the origin,  for which the multiset  $\overline\cS_m$ of allowed steps has generating polynomial
  \[
    \overline S_m(u,v):= \sum_{(a,b) \in \overline\cS_m}u^a v^b= \frac{(1+u)^m(1+v)(1+u+v)}{uv}.
  \]
\end{corollary}

\begin{proof}
  Let us introduce a new trivariate series $R_m(t;u,v)\equiv R_m(u,v)$ defined by $R_m(u,v):=Q'_m(1+u,1+v)$. Observe in particular that $R_m(0,0)= Q'_m(1,1)=G'_m(1,1)/t$ is  the \gf\ of  intervals in the posets $\LD'_{m,n}$, where $t$ records the size of the paths, minus $1$. The equation of Proposition~\ref{prop:eq-down} gives:
  \[
    R_m(u,v)=1+t \frac{(1+u)^m(1+v)(1+u+v)}{uv} R_m(u,v) - t \frac{(1+v)^2}{uv} R_m(0,v)-t \frac{(1+u)^{m+1}}{uv} R_m(u,0)+t \frac {R(0,0)} {uv} .
  \]
  This equation precisely describes (see~\cite{BoBoMe18}) the \gf\ of quadrant walks with steps in $\overline\cS_m$ starting from the origin, where $t$ records the length (number of steps), and $(u,v)$ the final coordinates.
\end{proof}

In particular, for $m=1$ the multiset  of allowed steps has generating polynomial 
    \[
      \overline S_1(u,v)=3+ \frac 1{uv} + \frac 2 u + \frac 2 v + \frac u v + \frac v u +u+v,
    \]
    and, up to the (harmless) constant term $3$, coincides with the steps of an already solved quadrant model~\cite[Sec.~A.7]{BeBMRa-17}. In the next section, we take advantage of this earlier result to prove Theorem~\ref{thm:counting}. We also explain how the principles of~\cite{BeBMRa-17} can be used to solve, instead of the equation obtained for $R_1(u,v)$ (or equivalently, for $G'_1(x,y)$), the equation obtained in Proposition~\ref{prop:eq-up} for $G_1(x,y)$.

\section{Exact enumeration of   intervals in $\LD_n$}
\label{sec:m=1}
In this section we focus on the case $m=1$, and drop all indices $m$ in our \gfs, writing for instance $G(x,y)$ instead of $G_1(x,y)$. Our aim is to prove, and actually refine, Theorem~\ref{thm:counting} that gives the size \gf\ of intervals in the lattices $\LD_n$. We can choose between two starting points to address this question. Either we start from the equation for $Q(x,y)=G(x/y,y)/(tx)$ in Proposition~\ref{prop:eq-up}:
 \beq\label{eq:Q1}
    Q \! \left(x , y\right)=1  +tx Q \! \left(x , y\right)
    +t  y^2\, \frac{x Q \! \left(x , y\right)-yQ \! \left(y, y\right)}{\left(x -y\right) \left(y -1\right)}
    -t  \,\frac{x Q \! \left(x  , 1\right)-Q \! \left(1, 1\right)}{\left(x -1\right) \left(y  -1\right)},
  \eeq
or from the equation for $Q'(x,y)= G'(x,y)/(txy)$ in Proposition~\ref{prop:eq-down}:
 \[
    Q'(x,y)=1+tx\frac{yQ'\! \left(x , y\right)-Q' \! \left(x, 1\right)}{y-1}
    +t y^{2}\, \frac{ x Q' \! \left(x , y\right)-Q' \! \left(1, y\right)}{\left(x -1\right) \left(y -1\right)}
    -t\, \frac{ x Q' \! \left(x , 1\right)-Q' \! \left(1, 1\right)}{\left(x -1\right)\left(y -1\right) }
    .\]
Recall from Corollary~\ref{cor:finite} that up to the change of variables $(x,y) \mapsto (1+x,1+y)$, the latter equation is equivalent to 
 \beq\label{eq:R1}
    R(x,y)=1+t \frac{(1+x)(1+y)(1+x+y)}{xy} R(x,y) - t \frac{(1+y)^2}{xy} R(0,y)-t \frac{(1+x)^{2}}{xy} R(x,0)+t \frac {R(0,0)} {xy} .
  \eeq
Again,   note that this equation is symmetric in $x$ and $y$, and hence the same holds for the series $R(x,y)$ and $G'(x,y)$.

In all three cases (and in fact also for $m>1$) we  have a linear equation, relating a main unknown series $F(x,y)$ to some specializations involving at most one of the two \emm catalytic, variables~$x$ and $y$: for instance $F(x,1)$, $F(1,1)$, $F(y,y)$ or $F(x,0)$. Historically, the first (non-linear) equations of this type appeared in the seventies in Tutte's work on the enumeration of properly coloured planar maps~\cite{lambda12,tutte-dichromatic-sums}. He devoted ten years and ten papers to the solution of just one of them (see~\cite{tutte-chromatic-revisited} and references therein). More recently, similar equations appeared in more contexts, like the enumeration of certain classes of permutations~\cite{bousquet-motifs,bouvel-baxter}, or of lattice walks confined to a quadrant~\cite{BeBMRa-17,BoKa08,BoRaSa14,bomi10,DHRS-17,raschel-unified}, or of some classes of maps~\cite{albenque-menard-schaeffer,BeBM-11,BeBM-15,BMRF-bipolar}. The solutions of these equations are not systematically algebraic (nor even D-finite as will be proved here for $m>1$), but our series $Q$, $Q'$ and $R$ \emm will, be proved to be algebraic. Several, often \emm ad hoc, approaches have been designed to prove algebraicity for such equations~\cite{BoKa08,bous05,Bousquet2016elementary,bousquet-motifs,bousquet-versailles,bomi10,Mishna-jcta}. The most systematic one, based on a certain notion of \emm invariants,, is based on ideas developed by Tutte in his enumeration of coloured maps.
In particular, the invariant approach has already solved in~\cite{BeBMRa-17} an equation that is extremely close to the above equation defining $R(x,y)$, and we will rely on this to solve~\eqref{eq:R1} and determine $G'(x,y)$ (Proposition~\ref{prop:Q1p-complete}). Next, we will see that an invariant approach can also be applied to determine $Q(x,y)$ (Proposition~\ref{prop:Q1-complete}), even if Eq.~\eqref{eq:Q1}, due to the terms $Q(y,y)$ and the denominator $(x-y)$, looks  rather different from~\eqref{eq:R1}.

  \subsection{The first ascent of $\bm P$, and the statistics $\bm{r(P,Q)}$}
\begin{proposition}\label{prop:Q1p-complete}
  The series $G'(t;x,y)$ that counts intervals $[P,Q]$ in the lattices $\LD_n$ by the size,
  the height of the first ascent of $P$
  and the statistics $r(P,Q)$ defined above Lemma~\ref{lem:first-insert} 
  is symmetric in~$x$ and~$y$. It is algebraic of degree $12$ over $\qs(t,x,y)$, and can be expressed as follows.

  Let $Z$  be the only series in~$t$ with constant term $0$ satisfying
  \[
Z=t(1+Z)(1+2Z)^2.
  \]
  Then the size \gf\ of intervals in the ascent lattices $\LD_n$ is the following cubic series:
  \beq\label{G11Z}
    G'(1,1)= Z(1-2Z+2Z^3).
  \eeq
  More generally the bivariate series $G'(x,1)$ is
\[
   G'(x,1)=G'(1,x) =\frac {C_0(x) - C_1(x) \sqrt{\Delta(x)}}{2x^3(1-x)Z^2},
 \]
 where
  \begin{align*}
    \Delta(x)&=  \left(1+Z \right)^{2} \left(1+2 Z \right)^{2}
    -2 Z \left(Z +1\right) \left(2 Z^{2}+4 Z +1\right) x
    +Z^{2} x^{2},
  \\
  C_1(x)&=\left(\left(1+2 Z \right)^{2}-2 Z^{2} x  -Z \,x^{2}\right)
  \left(\left(1+Z \right) \left(1+2 Z \right)-2 \left(1+Z \right)^{2} x +x^{2}\right),
  \end{align*}
   and
  {\small \begin{multline*}
      C_0(x)=\left(1+Z \right)^{2} \left(1+2 Z \right)^{4}
      -\left(1+Z \right) \left(1+2 Z \right)^{2}\left(8 Z^{3}+16 Z^{2}+9 Z +2\right)  x
      -Z \left(1+Z \right) \left(1+2 Z \right)^{2} x^{4}     \\
       +Z^{2} x^{5}+\left(1+Z \right) \left(18 Z^{5}+46 Z^{4}+48 Z^{3}+25 Z^{2}+7 Z +1\right) x^{2}
      -Z \left(2 Z^{5}-6 Z^{3}-7 Z^{2}-3 Z -1\right) x^{3}.
  \end{multline*}}
An algebraic expression of $G'(x,y)$ in terms of $x,y$ and $Z$ can then be obtained from the functional equation of Proposition~\ref{prop:eq-down}.
 \end{proposition}

 \noindent{\bf Remarks}\\
 {\bf 1.} The above expression of $G'(1,1)$ coincides with the value of $\sum_ng(n) t^n$ given in Theorem~\ref{thm:counting}. Deriving from this the asymptotic behaviour of the numbers $g(n)$ is a routine task, following the principles of \emm singularity analysis,~\cite[Sec.~VII.7]{flajolet-sedgewick}. The series $Z$ is found to have a unique dominant singularity, of the square root type, located  at $t_c:=1/\mu$ where $\mu$ is given in Theorem~\ref{thm:counting}.
 One needs to expand~$Z$ around $t_c$ up to the order of $(1-\mu t)^{5/2}$ to obtain from~\eqref{G11Z} the singular behaviour of $G'(1,1)$, in the form
 \[
   G'(1,1)= c_0 +c_1 (1-\mu t) +c_2(1-\mu t)^2 +c_{5/2} (1-\mu t)^{5/2} + \LandauO\left((1-\mu t)^3\right).
\]
The constant $\kappa$ of Theorem~\ref{thm:counting} is then $c_{5/2}/\Gamma(-5/2)$.
\\
{\bf 2.} Since any algebraic series is also D-finite, the series $G'(1,1)=\sum_n g(n) t^n$ also satisfies a linear differential equation. The corresponding linear recurrence relation 
reads
\[
  \left(n +4\right) \left(2 n +7\right) g \! \left(n +2\right)=
2 \left(11 n^{2}+44 n +42\right) g \! \left(n +1\right)+  n \left(2 n +1\right) g \! \left(n \right).
\]
One may wonder whether this can be explained combinatorially directly on intervals.\\
{\bf 3.} In the same way we have parametrized rationally $t$ and $G'(1,1)$ by $Z$, we can write $x$ and $G'(x,1)$ as rational functions in $Z$ and the unique series $U(x)$ with constant term $0$ satisfying
\[
  U(x)=   t x(1+U(x) )\left(1+3Z+Z^2 + Z(1+Z) U(x) \right).
\]
Note that $U(1)=Z/(1+Z)$. One can then write  $x$ rationally  in terms of $U$ and $Z$. In particular, the discriminant $\Delta(x)$ becomes a square, and finally
\[
  G'(x,1)=G'(1,x)=\frac{tx(1+U)}{(1+2Z)^{2} \left(1+3Z+Z^2-Z^2U\right)}\  P(Z,U)
  \]
 with {\small \begin{multline*}
    P(Z,U)=-Z^{7} \left(1+Z\right)^{3} U^{4}-Z^{5} (1+Z)^2\left(2 Z^{3}+3 Z^{2}-4 Z -2\right) U^{3}\\
     +Z^3(1+Z)\left(13 Z^{5}+49 Z^{4}+62 Z^{3}+24 Z^{2}-1\right) U^{2}\\
     +Z^2 \left(2 Z^{8}+35 Z^{7}+140 Z^{6}+221 Z^{5}+120 Z^{4}-48 Z^{3}-80 Z^{2}-31 Z -4\right) U \\
     +\left(Z^2+3Z +1\right) \left(Z^{8}+15 Z^{7}+31 Z^{6}+10 Z^{5}-19 Z^{4}-7 Z^{3}+9 Z^{2}+6 Z +1\right).
   \end{multline*}}

\begin{proof}[Proof of Proposition~\ref{prop:Q1p-complete}]
  As discussed at the end of Section~\ref{sec:recursive}, the series $R(t;x,y)=G'(t;1+x,1+y)/(t(1+x)(1+y))$ counts, by length and final coordinates, quadrant walks with step polynomial
     \[
      \overline S(x,y)=3+ \frac 1{xy} + \frac 2 x + \frac 2 y + \frac x y + \frac y x +x+y.
    \]
    These walks have already been counted in~\cite{BeBMRa-17}, but without the empty step, which has multiplicity~$3$ here. This means that the series $R$ is related to the series denoted $Q$ in~\cite{BeBMRa-17} and $\Qo$ here, to avoid confusion, by
    \[
      R(t;x,y)=\frac 1 {1-3t} \Qo\left( \frac t{1-3t}; x,y\right).
    \]
    We now dig in the details of the solution presented in~\cite[App.~A.7]{BeBMRa-17}. The series denoted~$Z(t)$ in Eq.~(A.8) of~\cite{BeBMRa-17}, and $\Zo(t)$ in the present paper, is related to the series $Z(t)$ of the proposition by
    \[
      Z(t)= 2 \times  \Zo\left(\frac t{1-3t}\right).
    \]
    Using Eq.~(A.9) of~\cite{BeBMRa-17}, which gives the expression of $\Qo(t;0,0)$, we then  obtain
    \[
      G'(t;1,1)=t R(t;0,0) = \frac 1 {1-3t} \Qo\left( \frac t{1-3t}; 0,0\right)= Z(1-2Z+2Z^3).
    \]
    Then Eq.~(A.10) in~\cite{BeBMRa-17} gives a quadratic equation for $\Qo(t;0,y)$, or equivalently $\Qo(t;y,0)$, with coefficients that are rational expressions in $t$, $y$ and two series $A_1$ and $A_2$, having themselves rational expressions in $\Zo$. From this we derive a quadratic equation for $G'(t;x,1)$ in terms of~$x$ and $Z$ (details are given in our {\sc Maple} session). Solving this equation gives the announced rational expression of $G'(t;x,1)$ in terms of $x$, $Z$ and $\sqrt{\Delta(x)}$.
\end{proof}

\subsection{The finals descents of $\bm P$ and $\bm Q$}

  We will now solve the equation~\eqref{eq:Q1} defining $Q(x,y)$,  and thus count intervals $[P,Q]$ by the heights of the final descents of $P$ and $Q$. Our solution is inspired from the \emm invariant approach, used in~\cite{BeBMRa-17} to solve, among other equations, the equation~\eqref{eq:R1} defining $R(x,y)$ (or more precisely, its variant with no empty step). The proof is detailed in Appendix~\ref{app:G11}. The specialists of equations in two catalytic variables may be interested in the fact that it uses both an \emm additive, and a \emm multiplicative decoupling,.

 \begin{proposition}\label{prop:Q1-complete}
 The series $G(t;x,y)$ that counts intervals $[P,Q]$ in the lattices $\LD_n$ by the size,
  the height of the last descent of $P$ and  the height of the last descent of $Q$,
 is algebraic of degree $12$ over $\qs(t,x,y)$, and can be expressed as follows.

  Let $Z$  be the only series in~$t$ with constant term $0$ satisfying
  \[
Z=t(1+Z)(1+2Z)^2.
\]
Then
\[
  G(x,1)= \frac{C_0(x)-C_1(x) \sqrt{\Delta_1(x)}}{2x^2Z^2},
\]
with
\begin{align*}
  \Delta_1(x)&=(1+Z-xZ)\left( \left(1+2 Z \right)^{2} \left(1+Z \right)-x Z\right),  \\
  C_1{(x)}&= (x-1)\left(\left(1+2 Z \right)^{2}-2 x Z \right)
             \left(2 \left(1+Z \right) \left(1+2 Z \right)-x\right)
\end{align*}
and
\allowdisplaybreaks
\begin{multline*}
  C_0{(x)}=
-2\left(1+Z \right)^{2} \left(1+2 Z \right)^{4}   +3 \left(1+Z \right) \left(1+2 Z \right)^{2}\left(4 Z^{3}+8 Z^{2}+6 Z +1\right)  x -
  \\
  \left(12 Z^{6}+64 Z^{5}+132 Z^{4}+134 Z^{3}+70 Z^{2}+16 Z +1\right) x^{2}
 +Z \left(12 Z^{3}+24 Z^{2}+16 Z +3\right) x^{3} -2 Z^{2} x^{4}.
\end{multline*}
Analogously,
\[
  G(1,y)= \frac{D_0{(y)}-D_1{(y)} \sqrt{\Delta_2(y)}}{2yZ(1-y)^2},
\]
with
\begin{align*}
  \Delta_2(y)&= (1 + 2Z)^2-4yZ(1 + Z),  \\
  D_1{(y)}&= (1+2Z-y)(2+2Z-y(1+2Z)),
\end{align*}
and
\[
  D_0{(y)}=-2 Z \,y^{3}+\left(4 Z^{3}+8 Z^{2}+10 Z +1\right) y^{2}
  -\left(12 Z^{3}+24 Z^{2}+16 Z +3\right) y +2 \left(1+2 Z \right)^{2} \left(1+Z \right).
\]
An algebraic expression of $G(x,y)$ in terms of $x,y$ and $Z$ can then be obtained from the functional equation of Proposition~\ref{prop:eq-up}.
\end{proposition}

\noindent {\bf Remarks}\\
{\bf 1.} Note that the second factor in $\Delta_1(x)$ is simply $Z(1-tx)/t$.\\
{\bf 2.} Again we have rational parametrisations for the series $G(x,1)$ and $G(1,y)$, by series ${U\equiv}U(x)$ and ${V\equiv}V(y)$ (with no constant term in $t$) satisfying
\[
U= \frac{xZ(1+U)}{(1+Z)(1+2Z-Z^2U)}, \quad \text{and} \quad  V=  \frac{yZ(1+Z)(1+V)^2}{(1+2Z)^2},
\]
respectively. 
Indeed,
\[
  G(x,1)= U\frac{(1+Z) \,P(Z,U)}{(1+2Z-Z^2U)^2},
\]
with
\begin{multline*}
  P(Z,U)=2 Z^{6} \left(1+Z \right) U^{3}+Z^{4} \left(8 Z^{3}+4 Z^{2}-8 Z -5\right) U^{2}\\+Z^{2} \left(6 Z^{5}-6 Z^{4}-20 Z^{3}-5 Z^{2}+9 Z +4\right) U -\left(1+2 Z \right) \left(4 Z^{5}+4 Z^{4}-2 Z^{3}-3 Z^{2}+Z +1\right),
 \end{multline*}
and
\[
  G(1,y)=V \, \frac {2 Z^{2} \left(1+Z \right)^{2} V^{2}+Z \left(4 Z^{3}+4 Z^{2}-4 Z -3\right) V +\left(1+Z \right) \left(2 Z^{3}-2 Z^{2}+1\right)}{(1+Z-ZV)^2}.
\]

\section{Asymptotic enumeration of  intervals}
\label{sec:asympt}
Our aim in this section is to prove Propositions~\ref{prop:Dmn} and~\ref{prop:Dmnp}, which give asymptotic estimates for the interval numbers in $\LD_{m,n}$ and $\LD'_{m,n}$ and establish non-D-finiteness results.
The key ingredient is the existence of bijections between intervals and quadrant walks, described in Section~\ref{sec:recursive}, combined with general asymptotic results on such walks by Denisov and Wachtel~\cite{denisov-wachtel} and their application to enumeration by Bostan et al.~\cite{BoRaSa14}. So far these results do not seem to have been applied to walks with infinitely many allowed steps, but, as will shall see by following the arguments of~\cite{BoRaSa14}, this does not raise difficulties as long as the step \gf\ converges in a sufficiently large domain.

\subsection{Intervals of $\LD_{m,n}$}
\begin{proof}[Proof of Proposition~\ref{prop:Dmn}]
  According to the second part of Corollary~\ref{cor:bij-mn}, the number $g_m(n)$ of intervals in the ascent lattice $\LD_{m,n}$ is also the number of quadrant excursions of length $n$ taking their steps in the set $\cS_m$. The \gf\ of this set is
  \[
    S_m(x,y)= x^m+ \sum_{i=-\infty}^{m-1} x^i \sum_{j=-\infty}^{m-i} y^j
    = x^m+ \sum_{i=-\infty}^{m-1} x^i \frac{y^{m-i}}{1-1/y}=\frac{x^{m} \left(x y -x +y \right)}{\left(x -y \right) \left(y -1\right)}.
  \]
  It converges absolutely for $1<|y|<|x|$. We now follow the probabilistic arguments of~\cite[Sec.~2.3]{BoRaSa14} (see also~\cite[Sec.~1.5]{denisov-wachtel}). We consider a random walk $(Y_1(n),Y_2(n))_{n\ge 0}$ in $\zs^2$, starting at $(0,0)$ and taking its steps in $\cS_m$, where each step $(i,j)$ occurs with a probability  ${x_0}^i{y_0}^j/S_m(x_0,y_0)$, for some $(x_0,y_0)$ chosen such that $1<y_0<x_0$. We moreover require that
  \[
    \frac {\partial S_m}{\partial x}(x_0,y_0)=    \frac {\partial S_m}{\partial y}(x_0,y_0)=0,
  \]
  as this choice guarantees that the walk $(Y_1,Y_2)$ has no drift, that is, the average displacement is zero. This gives
  \[
    x_0=\frac{2+\sqrt {m^2+4}}m, \qquad y_0=\frac{\sqrt{m^2+4}-m+2}2= \frac m 2 (x_0-1),
  \]
  and one can check that indeed, $1<y_0<x_0$.  The next step is to apply to the walk $Y=(Y_1,Y_2)$ a linear transformation so that the covariance matrix of the resulting walk, denoted $Z=(Z_1, Z_2)$, is the identity. As in~\cite{BoRaSa14}, the image by this transformation
  of the first quadrant ends up being
  a wedge $W_c$ of opening $\arccos(-c)$, where
  \[
    c= \frac{\frac{\partial^2 S_m}{\partial x \partial y}}
    {\sqrt{\frac{\partial^2 S_m}{\partial x^2} \frac{\partial^2 S_m}{\partial y^2}  }}(x_0,y_0).
  \]
  This coincides with the value of $c$ given in Proposition~\ref{prop:Dmn}. Now the probability that the walk $(Z_1,Z_2)$ visits $(0,0)$ at time $n$ without exiting the wedge $W_c$ before is
  \[
 p_m(n):=   \frac{g_m(n)}{S_m(x_0,y_0)^n},
  \]
  as each of the corresponding trajectories has probability $1/S_m(x_0,y_0)^n$. Now the random walk~$Z$ satisfies the conditions of~\cite{denisov-wachtel}: its steps are not contained in any (linear) half-plane, it is aperiodic (since the step $(0,0)$ is allowed), has no drift, its covariance matrix is the identity, and there are finite moments of any order. Moreover, the point $(0,0)$ can be reached from infinity (an assumption that seems to be missing in~\cite{denisov-wachtel}, see~\cite[Sec.~3.3]{BoBoMe18}). By~\cite[Thm.~6]{denisov-wachtel}, there exists a positive constant $\kappa$ such that
  \[
    p_m(n)\sim \kappa n^{-1-\pi/\arccos(-c)}.
  \]
  Combining this with the previous identity, and using
  \[
    \mu:=S_m(x_0,y_0)=\frac{ m \sqrt{m^{2}+4}+m^{2}+2}{2}\cdot \left(\frac{2+\sqrt{m^{2}+4}}{m}\right)^{m},
  \]
  yields the announced  asymptotic estimate of $g_m(n)$.

  Let us now discuss the implications of this result on the nature of the \gf\ of these numbers.

  When $m=1$, we obtain $c=(1-\sqrt 5)/4$, so that the uncorrelated random walk $Z$ lives in a cone of opening $\arccos(-c)=2\pi/5$, and the exponent $\alpha$ is $-7/2$, as already established in the previous section.

  Let us now prove that $m=1$ is the only integer value of $m$ for which $\pi/\arccos(-c)$ is rational. By~\cite[Thm.~3]{BoRaSa14}, this implies that the series $\sum_n g_m(n)t^n=G_m(1,1)$ is not D-finite when $m>1$. If $\arccos(-c)$ was a rational multiple of $\pi$, say of the form $\pi-\theta$, then we would have $c=\cos \theta= (z+1/z)/2$ for $z=e^{i\theta}$ a root of unity. Since $c$ is a root of the polynomial
  \[
    P_m(u):=4\left(m^{2}+3\right) u^{4}-4\left( m^{2}+2\right) u^{2}+m^{2},
  \]
  the two solutions $z$ and $1/z$ of $c=(z+1/z)/2$ are roots of
  \[
    \overline P_m(v):=\left(m^{2}+3\right) v^{8}+4 v^{6}+2\left( m^{2}+1\right) v^{4}+4 v^{2}+m^{2}+3.
  \]
  So we only have to prove that this polynomial admits no root of unity --- that is, no cyclotomic factor --- for $m>1$.  There are exactly $18$ cyclotomic polynomials of degree at most $8$. The first one is of course $\phi_1(v)=v-1$, and the last one is $\phi_{30}(v)=v^8 + v^7 - v^5 - v^4 - v^3 + v + 1$. We then take each of these polynomials $\phi(v)$ one by one, and reduce $\overline P_m(v)$ modulo $\phi(v)$ to detect if it could be a multiple of $\phi(v)$ for some values of $m$. But the leading coefficient of the remainder (which is a polynomial in $m$) never vanishes when  $m$ is an integer larger than~$1$. For instance, for $\phi=\phi_{30}$, we find
  \[
    \overline P_m(v)\!\!\!\! \mod \phi_{30} =  -  \left(m^{2}+3\right) z^{7}+4 z^{6}+\left(m^{2}+3\right) z^{5}+\left(3 m^{2}+5\right) z^{4}+\left(m^{2}+3\right) z^{3}+4 z^{2}-\left(m^{2}+3\right) z,
  \]
  and the leading coefficient $-(m^{2}+3)$ has no integer root. We conclude that no root of $\overline P_m(v)$ is a root of unity, so that the exponent $\alpha$ is irrational, and the \gf\ of the numbers~$g_m(n)$ is not D-finite.
  \end{proof}

\subsection{Intervals of $\LD'_{m,n}$}
\begin{proof}[Proof of Proposition~\ref{prop:Dmnp}]
 We now argue in a similar fashion for mirrored $m$-Dyck paths.  According to the second part of Corollary~\ref{cor:bij-mnp}, the number $g'_m(n)$ of intervals in the ascent poset $\LD'_{m,n}$ is also the number of quadrant excursions of length $n$ taking their steps in the set $\cS'_m$. The \gf\ of this set is
  \[
    S'_m(x,y)= x^m\sum_{j=-\infty}^0 y^j
    + \sum_{i=-\infty}^{m-1} x^i \sum_{j=-\infty}^{1} y^j
    = \frac{x^m}{1-1/y} + \frac{x^{m-1}}{1-1/x}\cdot\frac{y}{1-1/y}
    =\frac{x^{m} y \left(x+y -1 \right)}{\left(x-1 \right) \left(y -1\right)}
.
  \]
  It converges absolutely for $1<|x|$ and $1<|y|$. The critical point $(x_0, y_0)$ is now 
  \[
    x_0=\frac{2 m^{2}+1+\sqrt{4 m^{2}+1}}{2 m^{2}}, \qquad y_0= \frac{2 m +1+\sqrt{4 m^{2}+1}}{2 m}= 1+m(x_0-1),
  \]
and one can check that indeed, $1<x_0$ and $1<x_0$. The opening angle is found to be  $\arccos(-c)$, where $c$ takes the value of Proposition~\ref{prop:Dmnp}. One obtains the asymptotic estimate of $g'_m(n)$ as before, using
  \[
    \mu:=S'_m(x_0,y_0)  =\left(2m+\sqrt{1+4m^2}\right) \left( \frac{1+\sqrt{1+4m^2}}{2m}\right)^{2m}.
  \]

  It remains to prove that the exponent $\alpha$ is irrational for $m>1$. The argument is the same as in the previous subsection, with now
  \[
    \overline P_m(v)=\left(3 m^{2}+1\right) v^{8}+4 m^{2} v^{6}+2\left( m^{2}+1\right) v^{4}+4 m^{2} v^{2}+3 m^{2}+1.
  \]
  \end{proof}
  \begin{remark}
    One could alternatively establish the asymptotic estimate of $g'_m(n)$ via the bijection with weighted walks  mentioned in Corollary~\ref{cor:finite}. Recall that these walks  use  finitely many steps only.
  \end{remark}

  \section{Final comments}
  \label{sec:final}
Several questions that have been investigated on other Dyck posets may be asked for the ascent poset. For instance, what is the height of a random vertex of $P$ or $Q$ in an interval $[P,Q]$ (see~\cite{chapuy-tamari-asympt})? Can one write a $q$-analogue of our functional equations recording the size of the longest chain from $P$ to $Q$ (see~\cite{mbm-fusy-preville,mbm-chapoton})? 

Conversely, it seems to be the first time that an order induced on \emm mirrored, $m$-Dyck paths is examined. For the ascent orders as for the other Dyck orders (e.g., Tamari) this seems less natural than studying $m$-Dyck paths because mirrored paths do not form an upper ideal (while $m$-Dyck paths do). However, this is how we discovered here the connection with sylvester classes of $m$-parking functions, and were able to count these, at least asymptotically. Would interesting results arise from other orders? 

\begin{figure}[htb]
    \centering
    	\scalebox{0.9}{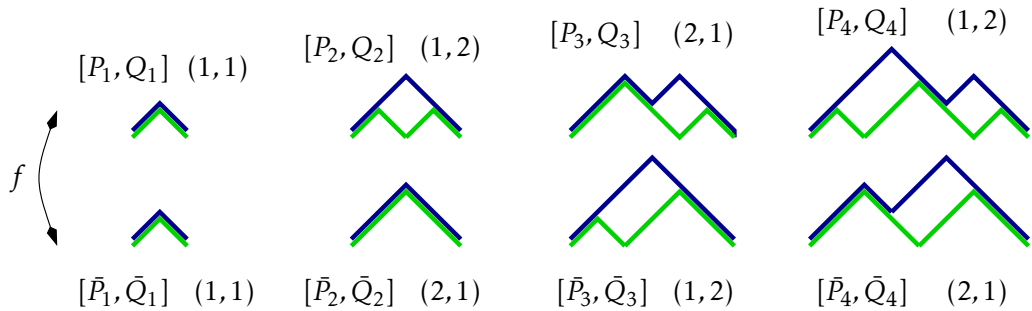}
    \caption{Recursive construction of an interval $[P_4, Q_4]$ in $\LD_4$, and of its image $[\bar P_4, \bar Q_4]$ by the involution $f$. For each interval $[P,Q]$ we give the statistics $(a(P), r(P,Q))$.}
    \label{fig:involution}
  \end{figure}

 The most immediate question raised by this paper is probably to explain bijectively, and non-recursively, the symmetry in $x$ and $y$ of the series $G'_1(x,y)$ that counts ascent intervals $[P,Q]$ in $\LD_1$ by the height $a(P)$ of the first ascent of $P$ (variable $x$) and the statistics $r(P,Q)$ defined at the beginning of Section~\ref{sec:recp} (variable $y$).  Recall that this is directly related to the fact that the quadrant walks that encode the recursive construction of these intervals (Corollary~\ref{cor:bij-mnp}) have a step set $\cS'_1$ that is $x/y$-symmetric. On these walks, the involution is obvious and consists in a reflection in the first diagonal. In recursive terms, if $w'$ is a quadrant walk and $w$ is the walk obtained by deleting its final step $(\delta_x,\delta_{y})$, then the image of $w'$ by the involution is obtained by appending the step $(\delta_{y}, \delta_{x})$ at the end of the image of~$w$. Recall the quadrant walk associated with an interval $[P,Q]$ such  $a(P)=a$ and $r(P,Q)=r$ ends at $(a-1, r-1)$. Hence we can define an involution $f$ on intervals $[P',Q']$, recursively on their size $n$ as follows:
\begin{itemize}
\item 
  $f([UD, UD])= [UD,UD]$
  \item  for $n>1$, let $b=a(P')$ and $s=r(P',Q')$, and let $[P,Q]$ be the interval of size ${n-1}$ obtained by deleting the initial peaks of $P'$ and $Q'$.  Then $f([P',Q'])$ is the only interval $[\bar P', \bar Q']$ obtained by insertion of peaks in  $f([P,Q])$ such that $a(\bar P')=s$ and $r(\bar P', \bar Q')=b$. An example is given is Figure~\ref{fig:involution}.
  \end{itemize}
  
\section*{Acknowledgements} We thank Clément Chenevière, Philippe Nadeau, and Benjamin Testart for interesting comments. \\
MBM was partially supported by the ANR projects DeRerumNatura (ANR-19-CE40-0018), Combin\'e (ANR-19-CE48-0011), and CartesEtPlus (ANR-23-CE48-0018). JLB and SK were supported in part by ANR PiCs (ANR-22-CE48-0002) and ANER ARTICO funded by Bourgogne-Franche-Comt\'e region.

\appendix

\section{Proof of Proposition~\ref{prop:Q1-complete}}
\label{app:G11}
Here we solve the equation~\eqref{eq:Q1} defining $Q(x,y)$, by adapting  the \emm invariant approach, of~\cite{BeBMRa-17}. All details of the calculations can be followed on an accompanying {\sc Maple} session,  available on the second author's webpage.

In Eq.~\eqref{eq:Q1}, let us group the terms involving $Q(x,y)$, and multiply out by $(y-1)$. We thus obtain
\beq\label{eq:H}
  K(x,y) (y-1) Q(x,y)= y-1  -t  \,\frac{xQ \! \left(x  , 1\right)-Q \! \left(1, 1\right)}{x -1}-\frac{ty^3}{x-y}  Q(y,y),
\eeq
where
\beq\label{ker:def}
  K(x,y)= 1- tx - \frac{txy^2}{(x-y)(y-1)}
\eeq
is  the \emm kernel, of the equation. 

\subsection{Invariants}

The notion of invariants is related to the expansion (in $t$) of $1/K(x,y)$. Namely,
\[
  \frac 1 {K(x,y)} = \sum_{n\ge 0}  t^n \left(x + \frac{xy^2}{(x-y)(y-1)}\right)^n.
\]
This is a series in $t$, with coefficients in $\qs(x,y)$. The denominators of these coefficients are powers of $(x-y)(y-1)$. Seen as fractions in $y$, they have a pole at $x$ and a pole at $1$, and the orders of these two poles increase with  the exponent of $t$. Let us now introduce a notion of series with \emm poles of bounded order,.

\begin{definition}
  Let $F(t;x,y)= \sum_n t^n f_n(x,y)$ be a Laurent series in $t$ with coefficients in $\qs(x,y)$.  We say that $F$ has \emm poles of bounded order at $y=1$ and $y=x$, if there exists an integer $m$ such that the coefficients of $(x-y)^m(y-1)^m F(t;x,y)$, seen as rational functions in $y$,  have no pole at $y=1$ nor $y=x$. We will often say, for short, that $F$ has \emm poles of bounded order,.
\end{definition}
Clearly, the above series $1/K(x,y)$  does \emm not, have poles of bounded order.

\begin{definition} A pair $(I(x),J(y))$  of Laurent series  in $t$ with coefficients in $\qs(x)$ and $\qs(y)$, respectively, is \emm a pair of invariants, (for the kernel $K(x,y)$) if the ratio
  \[
    \frac{I(x)-J(y)}{K(x,y)},
  \]
  expanded as a series in $t$ with coefficients in $\qs(x,y)$, has poles of bounded order.  
\end{definition}
Note that in this case, $J(y)$ itself has no pole at $y=x$ (because it does not depend on $x$), but may have a pole (of bounded order) at $y=1$. The following lemma allows us to build new pairs of invariants from old. Its proof mimics the proof of Lemma~2.8 in~\cite{mbm-tq-kreweras}.

\begin{lemma}\label{lem:linear} The componentwise sum of two pairs of invariants $\left(I_1(x), J_1(y)\right)$ and $\left(I_2(x), J_2(y)\right)$
  is still a pair of invariants. The same holds for their componentwise product.
\end{lemma}

The following lemma will be key to construct equations (for $Q(x,1)$, or $Q(y,y)$) in a \emm single, catalytic variable --- hence with an algebraic solution, by~\cite{BMJ06}.

\begin{lemma}\label{lem:inv}
  Let $(I(x),J(y))$ be a pair of invariants such that the ratio
    \[
    \frac{I(x)-J(y)}{(x-y)(y-1)K(x,y)},
  \]
  expanded in powers of $t$, has coefficients with no pole at $y=1$ nor $y=x$ (when seen as rational series in~$y$).
  Then $I(x)$ and $J(y)$ are equal, and in particular, they depend on $t$ only.
\end{lemma}
\begin{proof}
  Assume on the contrary that $I(x) \not = J(y)$, and write
  \[
    I(x)-J(y)=(x-y)(y-1) K(x,y) H(x,y),
  \]
  where $H(x,y)$ is a non-zero series in $t$, with coefficients in $\qs(x,y)$ and no pole at $y=1$ nor $y=x$. Let $m$ be the valuation of $H(x,y)$ in $t$, and denote by $h_m(x,y)\not = 0$ the  coefficient of $t^m$. Then, given that $K(x,y)=1+ \LandauO(t)$,
  \[
    I(x)-J(y)=(x-y)(y-1) h_m(x,y) t^m + \LandauO(t^{m+1}).
  \]
  Let us write $I(x)= \sum_n i_n(x) t^n$ and $J(y)=\sum_n j_n(y) t^n$.
 The above identity gives
  \[
    i_m(x)-j_m(y)= (x-y)(y-1) h_m(x,y).
  \]
  By assumption, $h_m(x,y)$, seen as a fraction in $y$, has no pole at $y=1$ nor at $y=x$. Hence the same holds for $j_m(y)$. Evaluating the above identity at $y=1$ shows that $i_m(x)=j_m(1)$, so that $i_m(x)$ does not depend on $x$. We now have
 \[
    j_m(1)-j_m(y)= (x-y)(y-1) h_m(x,y).
  \]
  Evaluating this at $y=x$ gives $j_m(x)=j_m(1)=i_m(x)$, but this forces $h_m(x,y)=0$, a contradiction.
\end{proof}

\subsection{A finite group, and rational invariants}
To the kernel $K$, given by~\eqref{ker:def}, one can associate as in~\cite{bomi10} a group of birational transformations of pairs $(u,v)$ that leave the value $K(u,v)$ unchanged. Solving $K(u,v)=K(u',v)$ for $u'$ gives $u'=u$ or $u'=v \left(uv -u +v \right)/(u-v)/(v-1)$. Solving $K(u,v)=K(u,v')$ for $v'$ gives $v'=v$ or $v'=uv/(uv-u+v)$. We introduce accordingly two transformations  $\Phi$ and $\Psi$ defined by
\[
\Phi(u,v)= \left(
\frac{v \left(uv -u +v \right)}{\left(u -v \right) \left(v -1\right)}
, v\right), \qquad
\Psi(u,v)= \left(u, \frac{u v}{uv -u +v}\right)  .
\]
One can check that they are involutions, and generate a group $\mathcal G$ of order $10$. The orbit of $(x,y)$ under the action of $\mathcal G$ is
\begin{multline*}
  (x, y){\overset{\Phi}{\longleftrightarrow}}
  \left(\frac{y \left(x y -x +y \right)}{\left(x -y \right) \left(y -1\right)}
    , y\right)
  {\overset{\Psi}{\longleftrightarrow}}
\left(\frac{y \left(x y -x +y \right)}{\left(x -y \right) \left(y -1\right)}
, \frac{x y -x +y}{x \left(y -1\right)}\right)
{\overset{\Phi}{\longleftrightarrow}}
\left(\frac{x y -x +y}{y \left(y -1\right)} ,
\frac{x y -x +y}{x \left(y -1\right)}\right)
{\overset{\Psi}{\longleftrightarrow}}\\
\left(\frac{x y -x +y}{y \left(y -1\right)}, \frac{x y -x +y}{y^{2}}\right)
{\overset{\Phi}{\longleftrightarrow}}
\left(\frac{x \left(x y -x +y \right)}{y \left(x -y \right)}
  ,\frac{x y -x +y}{y^{2}}\right)
{\overset{\Psi}{\longleftrightarrow}}
\left(\frac{x \left(x y -x +y \right)}{y \left(x -y \right)}, \frac{x}{x -y}
\right)   {\overset{\Phi}{\longleftrightarrow}}\\
\left(\frac{x y}{\left(y -1\right) \left(x -y \right)}, 
  \frac{x}{x -y}\right)   {\overset{\Psi}{\longleftrightarrow}}
\left(\frac{x y}{\left(y -1\right) \left(x -y \right)}, 
\frac{x y}{x y -x +y}\right)   {\overset{\Phi}{\longleftrightarrow}} \left(x, \frac{x y}{x y -x +y}\right)  {\overset{\Psi}{\longleftrightarrow}} (x,y).
\end{multline*}
 Groups of order $10$ also appear in the enumeration of several families of walks in the quadrant, including of course those used in the previous subsection~\cite{BeBMRa-17,KaYa-15}. The group $\mathcal G$ can be used to construct rational invariants. Guided by Theorem~4.6 in~\cite{BeBMRa-17}, we take any rational function $H(u,v)$, and compute the sum  of its values over all pairs $(u,v)$ of the above orbit. Denoting this sum by $H_\sigma(x,y)$, we now define $I_0(x)=H_\sigma(x,Y)$ and $J_0(y)=H_\sigma(X,y)$, where $X$ (resp. $Y$) is a root of $K(\cdot,y)$ (resp. $K(x, \cdot)$).  We could adapt the proof of Theorem~4.6 in~\cite{BeBMRa-17} to prove that $(I_0(x), J_0(y))$ then forms a pair of rational invariants. But this is not really needed, as we can simply apply the above recipe with some rational function $H(u,v)$, and check that the resulting pair is indeed a pair of invariants.

So let us start for instance with $H(u,v)=u$. The above recipe gives $I_0(x)=J_0(y)=-2+4/t$, which is a trivial pair of invariants. Starting instead from $H(u,v)=v$, we obtain a non-trivial pair, which, after dividing by $2$ and subtracting $2$, reads
\[
  I_0(x)=\frac 1 {1-tx} -\frac 1 {tx^2}+ \frac{1+t}{tx} +x(1-t)-tx^2,
\qquad   J_0(y)=-\frac{t}{(y-1)^2}+\frac{1-t}{y-1}- \frac 1 {ty^2} +\frac{1+t}{yt}+y.
\]
Indeed, one can check that
\beq\label{ratio0}
  \frac{I_0(x)-J_0(y)}{K(x,y)}= \frac{(x-y)(1-y+txy)(x+y-xy-xyt(1+x-xy))}{x^2y^2t(xt-1)(y-1)},
\eeq
which has poles of bounded order (at $y=1$ and $y=x$). Note that the pair $(I_0,J_0)$ does not satisfy the conditions of Lemma~\ref{lem:inv}.

\subsection{Decouplings, and a new pair of invariants}
Let us now return to the functional equation~\eqref{eq:H} that defines $Q(x,y)$. We would like to derive from it a new pair of invariants, that is, to transform it into an identity of the form
\beq\label{decoupled}
  K(x,y) H(x,y)= I(x)-J(y),
\eeq
where $H(x,y)$ has poles of bounded order. But there is an obstacle on our way, since  the term $ty^3/(x-y)Q(y,y)$  does not depend on $y$ only, but also on $x$. A similar  difficulty also arises, for instance, when counting quadrant walks weighted by interactions with the coordinates axes~\cite{beaton-owczarek-rechni}, or  walks avoiding a quadrant~\cite{mbm-tq-kreweras,dreyfus-trotignon,RaschelTrotignon2018Avoiding}, or in some continuous probabilistic models~\cite{BmEFHR-arxiv}.

  We are going to remedy this difficulty using a \emm multiplicative decoupling, (and later, an \emm additive decoupling,). First, we observe that by definition~\eqref{ker:def} of the kernel,
  \[
    \frac{ty^3}{x-y}= \frac{1-tx}x \cdot y(y-1)- \frac{y(y-1)}x K(x,y).
  \]
  This is what we call a \emm multiplicative decoupling,: we have written the problematic term $ty^3/(x-y)$ as the \emm product, of a series in $x$ and a series in $y$, modulo the kernel $K$. This would not be possible for any term of course. This allows us to rewrite~\eqref{eq:H} as follows:
  \beq\label{eq:H-bis}
    \frac{K(x,y)(y-1)}{1-tx} \big( xQ(x,y)-yQ(y,y)\big) = \frac{x(y-1)}{1-tx} -\frac{tx}{1-tx} \frac{xQ(x,1)-Q(1,1)}{x-1}-y(y-1)Q(y,y),
  \eeq
  and our  problem is solved, since the term in $Q(y,y)$ no longer involves $x$.

  However, we have created a second difficulty, and our equation still does not look like~\eqref{decoupled}: the constant term, which was formerly $(y-1)$ in~\eqref{eq:H}, now mixes $x$ and $y$. This new problem will be solved as well if we can find rational functions $A(x)$, $B(y)$ and $H(x,y)$ such that
  \beq\label{dec}
    \frac{x(y-1) }{1-tx}= A(x)+B(y) + K(x,y) H(x,y),
  \eeq
  where $H(x,y)$ does not contain a factor $K(x,y)$ in its denominator. This is what we call an \emm additive decoupling,. Again, an arbitrary rational function does not have, in general, an  additive decoupling.

  There are two ways to look for a solution $(A,B)$ of~\eqref{dec} (see~\cite[Sec.~4.2]{BeBMRa-17}  and our {\sc Maple} session for details). The first one is by guessing (say, the fraction $B(y)$) and checking. Let us explain how this works. Denoting by $Y_0$ and $Y_1$ the two roots of $K(x, \cdot)$, we observe that~\eqref{dec} implies that
  \[
    \frac{x}{1-tx}= \frac{B(Y_0)-B(Y_1)}{Y_0-Y_1}.
  \]
  We can then start from an Ansatz on the form of $B(y)$ (fixing the number of poles and their orders, but not their values), form the divided difference $(B(Y_0)-B(Y_1))/(Y_0-Y_1)$, write it as a fraction in $x$ and $t$ (since it is a symmetric function of the two roots $Y_0$ and $Y_1$), and solve the resulting identity for the coefficients occurring in the Ansatz. This approach readily gives a solution,
  \[
    B(y)=-\frac y t+ \frac 1 {t(y-1)} -\frac 1 {t^2y},
  \]
  from which we derive
  \[
    A(x)=\frac{2+x}t+\frac 1{t^2x} + \frac 1 {t(tx-1)}.
  \]
  We can now check that
  \beq\label{dec:plus}
    \frac{x(y-1)}{1-tx} = A(x)+B(y) + K(x,y) \frac{(x-y)(1-txy)}{xyt^2(1-tx)}.
  \eeq
  The second approach is constructive, and  consists in applying Theorem~4.11 in~\cite{BeBMRa-17} (even though our kernel is \emm not, of the same form as the kernels studied in~\cite{BeBMRa-17}). This gives an alternative solution, differing from the above solution $(A,B)$  by a pair of invariants. More precisely, this new solution reads:
  \[
    A(x) + \frac 9 {5t} I_0(x) -\frac 7 {5t^2}(1+3t), \qquad
    B(y) - \frac 9 {5t} J_0(y) +\frac 7 {5t^2}(1+3t).
  \]
  
  Let us now combine the functional equation~\eqref{eq:H-bis} with the decoupling relation~\eqref{dec:plus}: we obtain
  \beq\label{ratio1}
    \frac{(x-y)(y-1)K(x,y)}{1-tx}
    \left( \frac{xQ(x,y)-yQ(y,y)}{x-y} - \frac{1-txy}{t^2xy(y-1)}\right) = I_1(x)-J_1(y),
  \eeq
  with
  \[
    I_1(x)= \frac{2+x}t+\frac 1{t^2x} + \frac 1 {t(tx-1)}- \frac {tx}{1-tx} \frac{xQ(x,1)-Q(1,1)}{x-1},
    \quad
    J_1(y)= \frac y t- \frac 1 {t(y-1)} +\frac 1 {t^2y}+y(y-1)Q(y,y).
  \]
  That is, we have found a second pair $(I_1,J_1)$ of invariants, this time in terms of the unknown series $Q$. It does not satisfy the conditions of Lemma~\ref{lem:inv}.

\subsection{Equations in one catalytic variable}
  Our aim is now to combine polynomially our two pairs of invariants $(I_0,J_0)$ and $(I_1, J_1)$, thanks to Lemma~\ref{lem:linear}, to form a new pair $(I(x),J(y))$ satisfying the conditions of Lemma~\ref{lem:inv}. This will imply that $I(x)$ and $J(y)$ are in fact both equal to a series $C$ depending on $t$ only. We will thus obtain an equation in only one catalytic variable for each of the two series $Q(x,1)=G(x,1)/(tx)$ and $Q(y,y)=G(1,y)/t$, just by writing $I(x)=C=J(y)$.

  Let us examine the ratios
  \[
    \frac{I_0(x)-J_0(y)}{(x-y)(y-1) K(x,y)} \qquad \text{and} \qquad \frac{I_1(x)-J_1(y)}{(x-y)(y-1) K(x,y)},
  \]
  derived from~\eqref{ratio0} and~\eqref{ratio1}, respectively. The coefficient of $t^n$ in each of these two series has no pole at $y=x$. This is obvious for the first ratio, and for the second, this relies on the fact that  the divided difference
  \[
    \frac{xQ(x,y)-yQ(y,y)}{x-y}
  \]
  has polynomial coefficients in $x$ and $y$. However, the coefficient of $t^n$ in the first (resp. second) ratio has  a double (resp. simple) pole at $y=1$. Accordingly, we observe a double (resp. simple) pole at $y=1$ in $J_0(y)$ (resp. $J_1(y)$). More precisely, the singular expansions  at $y=1$ of these two series are respectively
  \[
 J_0(y)=   -\frac{t}{(y-1)^2}+\frac{1-t}{y-1} + \LandauO(1)\qquad \text{and} \qquad J_1(y)=-\frac 1 {t(y-1)}+ \frac{ 1+t} {t^2} + \LandauO(y-1).
  \]
  It is then natural to introduce the series
  \[
    J(y)= J_0(y) + t^3 J_1(y)^2 -t(1+3t) J_1(y),
  \]
  which has no pole at $y=1$. By Lemma~\ref{lem:linear}, the pair $(I(x),J(y))$ forms a pair of invariants, if we define $I(x)$ by
 \[
    I(x)= I_0(x) + t^3 I_1(x)^2 -t(1+3t) I_1(x).
  \]
  Using~\eqref{ratio0} and~\eqref{ratio1}, we can then check that this new pair of invariants satisfies the condition of Lemma~\ref{lem:inv}. This implies that $I(x)=J(y)=C$ for some series $C$.  The latter series is easily identified by expanding the expression of $J(y)$ at $y=1$. This gives
  \[
    I(x)=2-4t-2t^2 Q(1,1)=J(y).
  \]
  We can now replace $I(x)$ and $J(y)$ by their expressions in terms of $(xQ (x,1)-Q(1,1))/(x-1)$ and $Q(y,y)$, respectively. This gives two polynomial equations,
  \[
    \Pol_1(Q(x,1),Q(1,1),t,x)=0,  \qquad\text{and} \qquad \Pol_2(Q(y,y),Q(1,1),t,y)=0.
  \]
  Each of them is a  polynomial equation in one catalytic variable only, quadratic in the main series (that is,  $Q(x,1)$ or $Q(y,y)$).

  \subsection{Algebraicity}
  It remains to solve these two equations in one catalytic variable. We will work, say, with the equation for $Q(y,y)$, which is a bit lighter than the other. It reads:
    \beq\label{eqcat-y}
    \Pol_2(Q(y,y),Q(1,1),t,y)=0,
    \eeq
    where
  \[
  \Pol_2(q,q_1,y,t) =
  y^2 t^2 \left(y -1\right)^{2} {q}^{2}+\left(y \left(2 y^{2}-5 y +1\right) t -\left(y -1\right) \left(y -2\right)\right) q+2 tq_1   +\left(y -1\right) \left(y -2\right)=0.
\]
This is a very simple instance of a polynomial equation in one catalytic variable~\cite{BMJ06}: it is (only) quadratic in $Q(y,y)$, and involves a single additional unknown series depending on $t$ only, namely $Q(1,1)$. In this case, the machinery of~\cite{BMJ06} reduces to Brown's quadratic method~\cite{Brown65,goulden-jackson}. We first examine whether there exists formal power series $Y\equiv Y(t)$ such that the first derivative of the above polynomial $\Pol_2$, evaluated at $(Q(Y,Y),Q(1,1),t,Y)$, vanishes. That is,
\[
  (Y - 1)(Y - 2)= 2t^2Y^2(Y - 1)^2Q(Y,Y) +t Y(2Y^2 - 5Y + 1).
\]
The form of this equation shows that two such series exist, with constant terms~$1$ and~$2$, respectively. Then Brown's result tells us that each of these two series is a double root of the discriminant of $\Pol_2$ with respect to its first variable, evaluated at $(Q(1,1),t,y)$. That is, a double root of
\begin{multline*}
  -4 t \,y^{5}+\left(-8 Q_{11} \,t^{3}-7 t^{2}+22 t +1\right) y^{4}+\left(16 Q_{11} \,t^{3}+18 t^{2}-40 t -6\right) y^{3}\\
  +\left(-8 Q_{11} \,t^{3}-7 t^{2}+26 t +13\right) y^{2}-4\left( t +3\right) y +4,
\end{multline*}
 where we have written $Q(1,1)=Q_{11}$. Since this polynomial in $y$ has multiple roots, its discriminant vanishes. This gives for $Q_{11}$ the following cubic equation:
\[
  64  t^{6} Q_{11}^{3}+16 t^{3} \left(11 t^{2}-18 t -1\right) Q_{11}^{2}+\left(161 t^{4}-452 t^{3}+238 t^{2}-28 t +1\right) Q_{11} +49 t^{3}-167 t^{2}+25 t =1.
\]
One can now check that if we introduce the series $Z$ of Theorem~\ref{thm:counting}, the above equation factors and yields $Q_{1,1}=(1+Z)(1+2Z)^2(1-2Z+2Z^3)$, so that $G(1,1)=tQ(1,1)=Z(1-2Z+2Z^3)$.

Let us now return to~\eqref{eqcat-y}. We can express $t$ and $Q_{11}$ as fractions in $Z$, and we now have a quadratic equation for $Q(y,y)$ with coefficients in $\qs(Z,y)$. Solving it then gives the announced expression of $G(1,y)=tyQ(y,y)$ in Proposition~\ref{prop:Q1-complete}.

We now proceed similarly with  the other equation in one catalytic variable,  $\Pol_1(Q(x,1), \allowbreak Q(1,1),t,x)=0$, which is quadratic in $Q(x,1)$.  We replace  $t$ and $Q_{11}$ by their rational expressions in $Z$, solve the resulting equation, and obtain the announced expression of $G(x,1)=txQ(x,1)$. \qed

\end{document}